\documentclass[a4paper,12pt, reqno]{amsart}
\usepackage{amssymb,amsmath,latexsym}
\usepackage{graphics, epsfig, epstopdf}
\usepackage{color}
\usepackage{comment}
\usepackage{hyperref}
\allowdisplaybreaks

\def\1{{\bf 1}}

\def\nn{\nonumber}
\def\bee{\begin{equation}}
\def\eee{\end{equation}}
 \def\sB {{\mathcal B}}

\def\R {{\mathbb R}}

\newtheorem{thm}{Theorem}[section]
\newtheorem{lemma}[thm]{Lemma}
\newtheorem{defn}[thm]{Definition}
\newtheorem{prop}[thm]{Proposition}
\newtheorem{corollary}[thm]{Corollary}
\newtheorem{remark}[thm]{Remark}
\newtheorem{example}[thm]{Example}
\numberwithin{equation}{section}

\def\qed{{\hfill $\Box$ \bigskip}}

\def\AA{{\mathcal A}}

\def\BB{{\mathcal B}}

\def\LL{{\mathcal L}}
\def\DD{{\mathcal D}}
\def\FF{{\mathcal F}}

\def\R{{\mathbb R}}

\def\E{{\mathbb E}}

\def\P{{\mathbb P}}

\def\wh{\widehat}
\def\wt{\widetilde}
\def\pf{\noindent{\bf Proof.} }

\def\arccot{\mathrm{arccot}}
\usepackage{hyperref}
\begin{document}
\title[]{Positive self-similar Markov processes obtained by resurrection}
\author{ Panki Kim \quad Renming Song \quad and \quad Zoran Vondra\v{c}ek}
\thanks{P. Kim: This research is  supported by the National Research Foundation of Korea(NRF) grant funded by the Korea government(MSIP) (No. NRF-2021R1A4A1027378).}
\thanks{R. Song: Research supported in part by a grant from
the Simons Foundation (\#960480, Renming Song)}
\thanks{Z. Vondra\v{c}ek: Research supported in part by the Croatian Science Foundation under the project 4197.}

 \date{}

\begin{abstract}
In this paper we study positive self-similar Markov processes obtained by  
 (partially) resurrecting a strictly $\alpha$-stable process at its  
first exit time from $(0,\infty)$. 
We construct those processes by using the Lamperti transform.
We explain their long term behavior and give conditions for absorption at 0 
in finite time. In case the process is absorbed at 0 in finite time,
we give a necessary and sufficient condition for the existence of a recurrent extension.
The motivation to study resurrected processes comes from the fact that their jump kernels may explode at zero. 
We establish sharp two-sided jump kernel estimates for a large class of  resurrected stable processes.

\end{abstract}
\maketitle



\smallskip

\section{Introduction}\label{s:intro}
A $[0,\infty)$-valued standard Markov process (see \cite{BG68}) $X=(X_t, \P_x)$, $t\ge 0$, $x\ge0$, is called a 
\emph{positive self-similar Markov process} 
(pssMp) if there exists $\alpha>0$ such that for any $x>0$ and $c>0$, 
the law of $(cX_{c^{-\alpha}t}:\, t\ge 0)$ under $\P_x$ 
is equal to the law of $(X_{t}:\, t\ge 0)$ under $\P_{cx}$. 
One refers to $\alpha$ as the self-similarity index.
We will say that $X$ is a pssMp with the origin as a trap (or that $X$ is absorbed at the origin) if once $X$ hits the origin it stays there forever. 
Self-similar  processes were introduced by Lamperti in \cite{Lam72} where he 
established 
a one-to-one correspondence between 
pssMps up to the first exit time from $(0, \infty)$
and possibly killed L\'evy processes.
A detailed description of this correspondence, usually called the Lamperti transform, is given in Section \ref{s:prelim}.

A canonical example of a pssMp with origin as a trap 
is an $\alpha$-stable process in $\R$ absorbed at the origin upon exiting $(0,\infty)$. To be more precise, let $\eta=(\eta_t)_{t\ge 0}$ be a strictly $\alpha$-stable process in $\R$, $\alpha\in(0,2)$. Its L\'evy measure has a density 
\begin{equation}\label{e:stable-density}
\nu(x)=c_+\, x^{-1-\alpha}\1_{(x>0)}+c_-\, |x|^{-1-\alpha}\1_{(x<0)}, \quad x\in \R,
\end{equation}
where $c_+, c_-\ge 0$ and $c_+=c_-$ if $\alpha=1$. Let 
$\rho:=\P(\eta_1\ge 0)=\P(\eta_1>0)$
 be the positivity parameter, and set $\wh{\rho}:=1-\rho$. The process $\eta$ will be parameterized so that
\begin{equation}\label{e:c+c-}
c_+=\frac{\Gamma(\alpha+1)}{\Gamma(\alpha\rho)\Gamma(1-\alpha\rho)}\, ,\quad c_-=\frac{\Gamma(\alpha+1)}{\Gamma(\alpha\wh{\rho})\Gamma(1- \alpha\wh{\rho})}.
\end{equation}
Throughout the paper, we will exclude the cases of only one-sided jumps.  More precisely, the set of permissible parameters $(\alpha, \rho)$ is given by
$$
\{(\alpha,\rho): \alpha\in (0,1), \rho\in (0,1)\}\cup \{(\alpha, \rho): \alpha\in (1,2), \rho\in (1-1/\alpha, 1/\alpha)\} \cup \{(\alpha, \rho)=(1,1/2)\},
$$
cf.~\cite[p.399]{KPW14}.
We denote by $\P_x$, $x>0$, the law of $\eta$ starting at $x$.

Let $\tau=\tau_{(0,\infty)}:=\inf\{t>0: \eta_t\in (-\infty, 0]\}$ be the first exit time of $\eta$ from $(0,\infty)$. At time $\tau$, 
we send the process to $0$ where it stays forever, and thus arriving at the process $X^{\ast}_t:=\eta_t\1_{(t<\tau)}$, $t\ge 0$. 
The process $X^{\ast}=(X^{\ast}_t, \P_x)$ is a pssMp of index $\alpha$, cf.~\cite[Section 3.1]{CC06}. If $T_0:=\inf\{t>0: X^{\ast}_t=0\}$, then $T_0=\tau_{(0,\infty)}<\infty$ and $X^{\ast}_{T_0-}>0$ a.s. 
We denote by $\xi^{\ast}$ the L\'evy process associated with $X^{\ast}$ through the Lamperti transform.
  
In this paper, we introduce a large class of positive self-similar processes that can be obtained by 
 (partially)  resurrecting 
the strictly $\alpha$-stable process $\eta$ at the first exit time $\tau$.
If $z=\eta_{\tau}<0$ is the position where $\eta$ lands at the exit from $(0,\infty)$, we return the process into 
 $[0,\infty)$ 
according to a probability distribution $p(z, \cdot)$. 
If the process is returned to 0, it stays there forever. 
More precisely, to ensure self-similarity, 
we consider probability kernels 
$p:(-\infty, 0)\times \BB([ 0, \infty))\to [0, 1]$ satisfying the scaling condition
\begin{equation}\label{e:p-scaling-measure}
p(\lambda z, \lambda A)= p(z,A)\, \quad \text{for all }z<0,  A\in 
\BB( [0, \infty))
\text{ and }\lambda >0.
\end{equation}
All such kernels arise in the following way: Let $\phi$ be a probability measure on 
$\BB([ 0, \infty))$. Then
\begin{equation}\label{e:p-phi-measure}
p(z,A):=
\phi(|z|^{-1}{A}),\quad   \text{for all } z<0 \text{ and } A\in 
\BB( [  0, \infty))
\end{equation}
satisfies \eqref{e:p-scaling-measure}.   
Conversely, if $p(\cdot, \cdot)$ satisfies \eqref{e:p-scaling-measure} and if we set $\phi(A)=p( -1, A)$, then $p(\cdot, \cdot)$ is of the form \eqref{e:p-phi-measure}. We call $p(\cdot, \cdot)$ the \emph{return kernel}. 
Note that it follows from \eqref{e:p-scaling-measure} that
$$
\mathfrak{p}:=1-p(z, \{0\})
$$
is independent of $z<0$.

Let $j(x,z):=\nu(z-x)$ be the jump kernel of $\eta$. 
Set
\begin{equation}\label{e:int-kernel-0-measure} 
q_0(x,A)  
:=\int_{(-\infty,0)}j(x,z)p(z,A)\, dz, \quad x>0, A\in 
\BB( [0, \infty)),
\end{equation}
and note that 
$$
q_0(x, \{0\})=\int_{(-\infty,0)}j(x,z)p(z,\{0\})\, dz=(1-\mathfrak{p})\int_{-\infty}^0 c_- (x-z)^{-1-\alpha}\, dz =(1-\mathfrak{p}) \frac{c_-}{\alpha}x^{-\alpha}.
$$
We define a \emph{resurrection kernel} $q$ as the restriction of $q_0(x, \cdot)$ to $(0,\infty)$: 
\begin{equation}\label{e:int-kernel-measure}
q(x,A)
:=\int_{(-\infty,0)}j(x,z)p(z,A)\, dz, \quad x>0, A\in 
\BB((0, \infty)).
\end{equation}
The idea behind the kernel $q_0(x, \cdot)$ 
is that if $x=\eta_{\tau_-}$, then instead of sending $\eta$ to the origin at time $\tau$  
(thus obtaining the pssMp $X^{\ast}$), 
with probability $1-\mathfrak{p}$ we send $\eta$ to the origin,  and with probability $\mathfrak{p}=p(1, (0,\infty))$ we restart (or resurrect) it  according to the normalized kernel $q(x, \cdot)$. 
If $\mathfrak{p}\in (0,1)$ we call the process \emph{partially resurrected}, and when $\mathfrak{p}=1$ we say that it is fully resurrected (or just resurrected). 
By Lemma \ref{l:q-always-density}, $q(x,\cdot)$ is absolutely continuous with respect to the Lebesgue measure, and its density will be denoted by $q(x,y)$, $x,y>0$. 
By resurrecting  according to the normalized density we arrive at a pssMp $\overline{X}$ absorbed at the origin
with  jump kernel $J(x,y):=j(x,y)+q(x,y)$. The precise construction will be carried out in Section \ref{s:res-proc} by means of the Lamperti transform. 

Examples of pssMp that can be obtained by such resurrection include the path censored process from \cite{KPW14} 
(also called the trace process on $(0, \infty)$ of the stable process) 
in which case the return kernel is equal to the Poisson kernel of $\eta$ with respect to 
$(-\infty,0)$; the process from \cite{DROV17, Von21} 
 related to nonlocal problems with Neumann boundary condition 
in which case the return kernel is equal to $j(z, y)dy/\int^\infty_0j(z, u)du$; 
and the absolute value process $|\eta |=(|\eta_t|)_{t\ge 0}$ in which case $p(z,A)=\delta_{-z}(A)$, see Subsection \ref{ss:examples-return}. 
In all three examples above $\mathfrak{p}=1$ 
so we have full resurrection of $\eta$. 
If $p(z,A)=(1-  \mathfrak{p})\delta_{0}(A)+\mathfrak{p}  \delta_{-z}(A)$,
we recover the ricocheted process from \cite{Bud18, KPV21}.

One of the key questions that we will address in this paper is the behavior of $\overline{X}$ at its 
absorption time at 0.
In case $\mathfrak{p}<1$ it is clear from the Lamperti trichotomy (see Section \ref{s:prelim})  that the partially resurrected process $\overline{X}$ will be absorbed at 0 in finite time by a jump. In case $\mathfrak{p}=1$, the  absorption time 
may be infinite or finite,
and in the latter case it turns out that $\overline{X}$ is continuously absorbed at 0.
We answer the question of finite or infinite lifetime 
by studying the behavior of the L\'evy process 
$\overline{\xi}=(\overline{\xi}_t, {\bf P}_x)$, $t\ge 0, x\in \R$, associated to
$\overline{X}$ through the Lamperti transform. Let $\overline{\Psi}$ denote the characteristic exponent of $\overline{\xi}$. 
We will write ${\bf P}_0$ as ${\bf P}$ and denote expectation with respect to ${\bf P}$ by ${\bf E}$. 
When $\mathfrak{p}<1$, it is obvious that 
${\bf E}\overline{\xi}_1=-\infty$.
Our first main result is about the finiteness of ${\bf E}|\overline{\xi}_1|$ and provides an explicit expression of ${\bf E}\overline{\xi}_1$ 
when $\mathfrak{p}=1$.

\begin{thm}\label{t:derivative-at-zero} 
Suppose $\mathfrak{p}=1$. 
It holds that ${\bf E}|\overline{\xi}_1|<\infty$ if and only if
\begin{equation}\label{e:phi-int-log}
\int_{(0,{\infty})}|\log y|\phi(dy)<\infty.
\end{equation}
In this case
\begin{equation}\label{e:derivative-at-zero}
 {\bf E}[\overline{\xi}_1]=i\overline{\Psi}'(0)
=\Gamma(\alpha)\frac{\sin(\pi \alpha\widehat\rho)}{\pi}
\left(\pi\cot(\pi\alpha\widehat\rho)+
\int_{(0,{\infty})}(\log y)\phi(dy)\right).
\end{equation}
\end{thm}

It follows from \cite[Theorem 7.2]{Kyp14} and the Lamperti trichotomy (see Section \ref{s:prelim}) that (i) 
If ${\bf E}[\overline{\xi}_1] \ge 0$,
then $\limsup_{t\to \infty}\overline{\xi}_t=+\infty$, hence the absorption time of $\overline{X}$ is infinite, and 
 (ii) If ${\bf E}[\overline{\xi}_1] < 0$, 
then $\lim_{t\to \infty}\overline{\xi}_t=-\infty$, hence the 
absorption time of $\overline{X}$ is finite $\P_x$-a.s.~and  $\overline{X}$ is continuously absorbed at 0.
Therefore, since $\alpha\wh{\rho}\in (0,1)$, to deduce the long term behavior of $\overline{X}$ it suffices to determine the sign of the expression in 
the parenthesis in \eqref{e:derivative-at-zero}, see Corollaries \ref{c:s}--\ref{c:s2}.

In case when $\overline{X}$ is 
absorbed at zero in finite time,
 we give a definitive answer on the existence of its recurrent positive self-similar extensions.
\begin{thm}\label{t:R_e} 
Suppose that $\overline{X}$  is absorbed at zero in finite time. 
If
\begin{equation}\label{e:kappa_0}
\kappa_0:=\sup\{\kappa\in 
(0,{\infty}):\, \int_{(0,{\infty})}u^{\kappa}\phi(du)<\infty\}
\in (0, \infty],
\end{equation}
then (1) $\overline{X}$ has a positive self-similar  recurrent extension which leaves 0 continuously;    
and (2)  there exists $\kappa^\ast\in (0, \alpha)$ such that, for any $\beta \in (0, \kappa^\ast)$,  $\overline{X}$ has a positive self-similar recurrent extension which leaves 0 by a jump associated 
with an excursion measure of the form $c\beta x^{-1-\beta}dx, x>0$. 

Conversely, if \eqref{e:kappa_0} does not hold, then $\overline{X}$ has no positive self-similar recurrent extension.
\end{thm}

Positive self-similar Markov processes and their associated L\'evy processes
have been extensively 
studied in the last 15 years. We mention here the Lamperti stable processes \cite{CC06, CPP10}, hypergeometric processes \cite{KPR10, KP13}, 
double hypergeometric L\'evy processes \cite{KPV21}, 
and $\beta$-processes \cite{Ku10}. The interest in those families of processes was mostly motivated  by the Wiener-Hopf factorization.

Our motivation for studying pssMps comes from our research program on the potential theory of Markov processes with jump kernels degenerate at the boundary.
In \cite{KSV21, KSV22}, we introduced a large class of symmetric Markov processes in $\R^d_+$ with jump kernels decaying at the boundary and systematically studied their potential theory.
The trace process of a symmetric $\alpha$-stable process on $\R^d_+$ is degenerate in the sense that its jump kernel blows up at the boundary, see \cite{BGPR21}.
The same feature is true also for the process studied in \cite{DROV17, Von21}.
In \cite{KSV22b} we studied the potential theory of a large class of  symmetric Markov processes
in $\R^d_+$ with jump kernels blowing up at the boundary. 
The main examples of such processes are 
rotationally symmetric $\alpha$-stable processes resurrected upon exiting the upper-half-space.
The current 
paper concentrates on the one dimensional case, but we do not assume that the processes are symmetric.  
The jump kernel $J(x,y)=j(x,y)+q(x,y)$ of our resurrected process exhibits unusual and interesting behavior when $y\to 0$. 
Depending on the  return kernel $p(z,\cdot)$, $J(x, y)$ may tend to $\infty$ at various rates as $y\to 0$. Since $j(x,y)$ is bounded away from the diagonal, this explosion is due to the  resurrection kernel.

In order to state the  precise  result about behavior of the resurrection kernel, 
 we first need a definition.
\begin{defn}
	{\rm Let $g:(0,\infty) \to (0,\infty)$ and $\beta_1, \beta_2 \in \R$.
		\begin{enumerate}
			\item[(i)] 
			We say that $g$ satisfies the lower weak scaling condition at zero $L_1(\beta_1)$ (resp. at infinity $L^1(\beta_1)$) if there exists $c\in(0,1]$ such that
			$$ 
			\frac{g(R)}{g(r)} \geq c \left(\frac{R}{r}\right)^{\beta_1} \quad \text{for all} \quad r\leq R< 
			1\;(\text{resp.}\;
			1\le r\leq R).
			$$
			\item[(ii)] We say that $g$ satisfies the upper weak scaling condition at zero $U_1(\beta_2)$ (resp. at infinity  $U^1(\beta_2)$) if there exists $C\in [1, \infty)$ such that
$$ 
\frac{g(R)}{g(r)} \leq C \left(\frac{R}{r}\right)^{\beta_2} \quad \text{for all} \quad r\leq R< 1\;(\text{resp.}\; 1\le r\leq R).
$$
					\end{enumerate}
}
\end{defn}
Here is our third main result in which we assume that 
the restriction of the measure $\phi$ to $(0, \infty)$ is absolutely continuous with respect to the Lebesgue measure
and, with slight abuse of notation, denote its density also by $\phi$. 
The notation $a\asymp b$ means that $c\le b/a \le c^{-1}$ for some $c\in (0,1)$.

\begin{thm}\label{t:estimates-of-q}
Suppose that the density $\phi$ is strictly positive.
(1) 
If $x \le y\le 5x$, then
$$
q(x,y)\asymp q(y,x)\asymp x^{-1-\alpha}\asymp y^{-1-\alpha}.
$$

\noindent
(2) Suppose $\phi$ satisfies the lower weak scaling 
condition $L_1(\beta_1)$ at zero  with $\beta_1>-1-\alpha $. 
Then for $5x\le y$, 
$$
q(y,x)\asymp (y-x)^{-1-\alpha}\int_{\frac{x}{y-x}}^{1}  \phi(t)\frac{dt}{t}
\asymp y^{-1-\alpha}
\int_{\frac{x}{y-x}}^{1}  \phi(t)\frac{dt}{t}.
$$
Further, 
if  
$\phi$ also satisfies the upper weak scaling 
condition $U_1(\beta_2)$ at zero with $\beta_2<0$, 
for $5x\le y$, 
$$
q(y,x)\asymp (y-x)^{-1-\alpha} \phi\big( \frac{x}{y-x}\big)\asymp y^{-1-\alpha} \phi\big( \frac{x}{y}\big).$$
(3) Suppose $\phi$ satisfies the upper weak scaling 
condition  $U^1(\gamma_2)$  at infinity with $\gamma_2<0 $. 
Then for $5x\le y$, 
\begin{equation}\label{e:estimate-q-21}
q(x,y)
\asymp (y-x)^{-1-\alpha}\int_0^{ \frac{y-x}{x}} t^{\alpha} \phi(t){dt} 
\asymp y^{-1-\alpha}\int_1^{\frac{y}{x}-1}t^{\alpha} \phi(t){dt}.
\end{equation}
Further, 
if  
$\phi$ also satisfies the lower weak scaling 
condition $L^1(\gamma_1)$ at infinity  with $\gamma_1<-1-\alpha$,
for $5x\le y$, 
$$
q(x,y)
\asymp (y-x)^{-1-\alpha} 
\big(\frac{x}{y-x}\big)^{-1-\alpha} \phi\big(\frac{y-x}{x}\big)
\asymp x^{-1-\alpha}\phi\big(\frac{y}{x}\big).
$$
\end{thm}

As a consequence of this theorem we can deduce, see Corollary \ref{c:estimates-of-q}, that the jump kernel $J(x,y)$ of the path censored process goes to $\infty$ at rate $y^{-\alpha\rho}$ as $y\to 0$, and that the jump kernel of the process with the resurrection kernel $j(z, y)dy/\int^\infty_0j(z, u)du$ goes to infinity at rate $\log(1/y)$ as $y\to 0$.

\medskip
Stable process conditioned to stay positive and censored stable process can be regarded as
resurrected stable processes, 
see \cite{Ber93} and  \cite[Remark 3.3]{KPW14}.
However, they do not fall
into the framework of resurrected stable processes of this paper. 
To cover these processes, we  introduce a larger class of pssMps in 
Section \ref{s:modified-jump} of this paper.
This larger class incudes stable processes conditioned to stay positive, stable processes conditioned to hit 0 continuously, 
and censored stable processes as examples.
The jump kernel $J(x,y)=j(x,y)+q(x,y)$ of our resurrected stable process can be regarded as a modification of the original kernel $j(x,y)$. 
The jump kernels of the class of pssMps in Section \ref{s:modified-jump} are of the more general form $j(x,y)\sB(x,y)$ where $\sB:(0,\infty)\times (0,\infty)\to (0,\infty)$. 
In case of multidimensional isotropic stable process, the  analogous procedure is quite standard when $\sB$ is bounded from below and above by two positive constants, leading to the so-called stable-like processes. The situation when the function $\sB(x,y)$ decays and vanishes at the boundary of the state space  was recently studied in \cite{KSV21, KSV22} in the multidimensional case of 
the upper-half-space in $\R^d$ and the symmetric jump kernel $j(x,y)=|x-y|^{-d-\alpha}$. 
The paper \cite{KSV22b} deals with the multidimensional case in the upper-half-space, the same jump kernel $|x-y|^{-d-\alpha}$, 
with $\sB(x,y)$ exploding at the boundary. The symmetric 1-dimensional case is also covered in that paper, 
but differently to the current paper which is mostly focused on self-similarity, the main concern of \cite{KSV22b} is on potential-theoretic questions.

{\bf Organization of the paper}:
In the next section we recall some preliminary results related to 
pssMps and their connection to L\'evy processes through the Lamperti transform,
in particular the Lamperti trichotomy and the relationship between infinitesimal generators of those processes. We also introduce two families of return kernels and show how some of examples for pssMp from literature fit into these families.

Section \ref{s:res-proc} is central to the paper. Starting from the general  return kernel described by the measure $\phi$, we first look at the regular step process with the  resurrection kernel $q(x,y)$ as its jump kernel, and its counterpart through the Lamperti transform -- a compound Poisson process $\chi$. We compute the L\'evy measure of $\chi$ and its characteristic function. For particular examples given by \eqref{e:phi-special} and \eqref{e:phi-special-exp} we obtain more precise expressions. Then we construct the resurrected process $\overline{X}$, 
and compute the characteristic exponent $\overline{\Psi}$  of the corresponding L\'evy process $\overline{\xi}$. 
In Subsection \ref{ss:X-lifetime} we study the behavior of $\overline{X}$ by analyzing the derivative of the characteristic exponent 
$\overline{\Psi}$  of $\overline{\xi}$ at zero and give a proof of Theorem \ref{t:derivative-at-zero}. In Subsection \ref{ss:res-proc} we prove Theorem \ref{t:R_e}. 
Finally, we provide several concrete examples illustrating 
the behavior of $\overline{X}$ at its absorption time.

In Section \ref{s:sym-int-kernel} we give a necessary and sufficient condition on $\phi$ making the  resurrection kernel symmetric,  i.e., $q(x,y)=q(y,x)$ for all $x,y>0$, cf.~Theorem \ref{t:q-symmetric}. If, in addition, the underlying  $\alpha$-stable process $\eta$ is symmetric, the resulting resurrected process $\overline{X}$ will be also symmetric. 

In Section \ref{s:estimates-q} we provide a proof of Theorem \ref{t:estimates-of-q}. The key estimates are given in Lemma \ref{l:estimates-of-q} where it is assumed that $|x-y|=1$. By using scaling of the resurrection kernel, this suffices to prove the theorem. A particular example of a slightly modified density from the family \eqref{e:phi-special} gives additional insight of possible 
behaviors at zero and at infinity, see Corollary \ref{c:estimates-of-q}.

In Section \ref{s:modified-jump} 
we put the resurrected process in a more general context of processes with modified jump kernel. Given a process with jump kernel $j(x,y)$, we multiply it by a function $\sB(x,y)$ which changes the behavior of the original kernel. 
We show that the jump kernel of the resurrected process can be thought of as being modified by the function which explodes at the boundary, namely at zero. We end the paper by looking at the symmetric $\alpha$-stable case modified by a function $\sB(x,y)$, 
and establish the behavior of the modified process at its lifetime.

For the reader's convenience, we summarize the known examples of resurrected stable processes  that can be included
in our framework. Path-censored stable processes, the processes studied 
in \cite{DROV17, Von21}, the absolute value of a stable process and the ricocheted 
stable processes are examples of our resurrected stable processes. 
Censored stable processes, stable processes conditioned to stay positive, 
and stable processes conditioned to hit 0 
continuously are included in the more general framework of Section \ref{s:modified-jump}.

{\bf Notation}:
We use ``:='' to indicate definitions. Define $a\land b := \min\{a, b\}$ and $a\vee b := \max\{a, b\}$.
We write $f\asymp g$ if $f,g$ are nonnegative functions, $c^{-1}g\le f\le cg$ with some constant $c \in (0,\infty)$. We call $c$ the   
comparability constant. Lower case letters 
$c_i, i=1,2,  \dots$ are used to denote the constants in the proofs
and the labeling of these constants starts anew in each proof.
 $\Gamma$ denotes the Gamma function defined as $\Gamma(x)=\int_0^\infty y^{x-1}e^{-y}d y$,
$\psi$ denotes the digamma function defined as $\psi(x)=\frac{d}{dx}\log \Gamma(x)$,
and $B$ denotes the beta function defined by $B(x,y)=\Gamma(x)\Gamma(y)/\Gamma(x+y)$.


\section{Preliminaries}\label{s:prelim}
\subsection{Lamperti correspondence}\label{ss:lam-corr}
We start this preliminary section by briefly describing the correspondence between positive self-similar Markov processes and 1-dimensional L\'evy processes, usually called the Lamperti transform. 
Let $\xi=(\xi_t, {\bf P}_x)$, $t\ge 0, x\in \R,$ be a possibly killed L\'evy process sent to $-\infty$ at death. Define the integrated exponential process $I=(I_t)_{t\ge0}$ by
$$
I_t:=\int_0^t e^{\alpha \xi_s} ds, \quad t\ge 0,
$$
and let $\varphi$ be its inverse:
$$
\varphi(t):=\inf\{s>0:\, I_s >t\}, \quad t\ge 0.
$$
For each $x>0$, define $\P_x:={\bf P}_{\log x}$ and 
$$
X_t:=\exp\{\xi_{\varphi(t)}\}\1_{(t < I_{\infty})}, \quad t\ge 0.
$$
Then $X=(X_t, \P_x)$, $t\ge 0, x>0$, is a pssMp  of index $\alpha$ with absorption time $\zeta=I_{\infty}$. 
(See \cite{PS18} for an  in-depth  analysis of the absorption time.)
Conversely, for a pssMp $X=(X_t, \P_x)$, $t\ge 0, x>0$, of index $\alpha$, let
$$
S_t:=\int^t_0X^{-\alpha}_udu
$$
and let $T(\cdot)$ be its inverse
$$
T(t):=\inf\{u>0:\, S_u >t\}, \quad t\ge 0.
$$
For any $x\in \R$, define ${\bf P}_x:=\P_{e^x}$ and $\xi_s:=\log X_{T_s}$. Then $\xi=(\xi_t, {\bf P}_x)$, $t\ge 0, x\in \R$,
is a possibly killed  L\'evy process.
Moreover, we have the following 
three exhausting scenarios (see \cite[Theorem 13.1]{Kyp14}) that we refer to as 
the {\it Lamperti  trichotomy}:
\begin{itemize}
	\item[(1)] 
	$\P_{x}(\zeta=\infty)=1$
	for all $x>0$ in which case $\limsup_{t\to \infty}\xi_t=\infty$; 
	\item[(2)] 
	$\P_{x}(\zeta<\infty, X_{\zeta-}=0)=1$
	for all $x>0$ in which case $\lim_{t\to \infty}\xi_t=-\infty$;
	\item[(3)] 
	$\P_{x}(\zeta<\infty, X_{\zeta-}>0)=1$
	for all $x>0$ in which case $\xi$ is killed at an independent exponentially distributed random time.
\end{itemize}
In case (2),
we will say that $X$ is continuously absorbed at $0$.
 
Now we recall a few facts about 1-dimensional L\'evy processes.
Let $\xi=(\xi_t, {\bf P}_x)$, $t\ge0, x\in \R$,  be a 1-dimensional L\'evy process  with  characteristic triple $(d,\sigma, \nu)$. 
We will write ${\bf P}_0$ as ${\bf P}$ and denote expectation with respect to ${\bf P}$ by ${\bf E}$. 
Then 
\begin{equation}\label{e:char-fn}
{\bf E} \left[e^{i\theta \xi_t}\right]=e^{-t\Psi(\theta)}, \quad \theta\in \R,
\end{equation}
where the characteristic exponent $\Psi$ is given by
\begin{equation}\label{e:char-exp}
\Psi(\theta)=d
i\theta +\frac{1}{2}\sigma^2 \theta^2 +\int_{\R}\left(1-e^{i\theta x}+i\theta x\1_{(|x|\le 1)}\right)\nu(dx), \quad \theta \in \R.
\end{equation}
Recall that $\xi_1$ has finite expectation if and only if  
$\int_{|y|\ge 1}|y|\nu(dy)<\infty$, cf.~\cite[Theorem 25.3, Example 25.12]{Sat14}. In this case,
by differentiating \eqref{e:char-fn} we get that ${\bf E}[\xi_1]=i\Psi'(0)$.

If $\xi$ is killed at an independent exponential time of parameter $q\ge 0$ (when $\xi$ is sent either to a cemetery $\partial$ or to $-\infty$), the characteristic exponent becomes $\wt{\Psi}(\theta)=\Psi(\theta)+q$. Thus for the killed L\'evy process $\wt{\xi}$, the killing rate is equal to $\wt{\Psi}(0)$. 

Let $\AA$ be the infinitesimal generator of the semigroup of $\xi$ (possibly killed at rate $q\ge 0$) acting on $C_0(\R)$ (continuous functions vanishing at infinity). Then, cf.~\cite[Theorem 31.5]{Sat14}, $C_0^2(\R)\subset \DD(\AA)$, and for $f\in C_0^2(\R)$,
$$
\AA f(x)=-qf(x) -
df'(x)+\frac{1}{2}\sigma^2 f''(x)+\int_{\R}\left(f(x+y)-f(x)-f'(x)y \1_{(|y|\le 1)}\right) \nu(dy). 
$$ 

Let $X$ be the pssMp of index $\alpha$ corresponding to the L\'evy process $\xi$.
Its infinitesimal generator $\LL$  can be described as follows (cf.~\cite[Theorem 1]{CC06}, where we take the usual 
cutoff function
$\ell(y)=y\1_{[-1,1]}(y)$): If $f:[0,\infty]\to \R$ is such that $f$, $xf'$ and $x^2 f''$ are continuous on $[0,\infty]$ then it belongs to the domain of $\LL$ and 
\begin{eqnarray}\label{e:pssMp-inf-gen}
\LL f(x)&=&-qx^{-\alpha}f(x) +x^{1-\alpha}\left(-
d +\frac{1}{2}\sigma^2\right) f'(x)+\frac{1}{2}\sigma^2 x^{2-\alpha}f''(x) \nonumber \\
& & +x^{-\alpha}
\int_{(0, \infty)}
\left(f(ux)-f(x)- x f'(x)(\log u)\1_{[-1,1]}(\log u)\right)
\mu(du),
\end{eqnarray}
where $\mu(du)=\nu(du)\circ \log u$. 
By the change of variables $y=\log u$, we get
\begin{eqnarray}\label{e:pssMp-inf-gen-2}
\LL f(x)&=&-qx^{-\alpha}f(x) +x^{1-\alpha}\left(-
d+\frac{1}{2}\sigma^2\right)f'(x)+\frac{1}{2}\sigma^2 x^{2-\alpha}f''(x) \nonumber \\
& & +x^{-\alpha}\int_{\R}\left(f(xe^y)-f(x)-xf'(x)y\1_{[-1,1]}(y)\right)\nu(dy),
\end{eqnarray}
which corresponds to the formula in \cite[p.~4]{PR13}.  
In case $\nu$ has a density (which we also denote by $\nu$), the integral in \eqref{e:pssMp-inf-gen-2} can be (after a change of variables) written in the form
$$
\int_0^{\infty}\left(f(z)-f(x)-xf'(x)(\log z/x)\1_{[-1,1]}(\log z/x )\right)\nu(\log z/x)\frac{dz}{z}
$$
showing that the intensity of jumps from $x$ to $z$ (i.e.~the jump kernel of $X$) is given by $z^{-1}\nu(\log z/x )$.

\subsection{Strictly stable process absorbed at 0 and censored process}\label{ss:sspcp}

Recall that $\eta=(\eta_t, \P_x)$ denotes a strictly $\alpha$-stable process in $\R$, $\alpha\in(0,2)$. Thus $\eta$ is a L\'evy process with characteristic exponent given by 
\eqref{e:char-exp},  where $\sigma=0$ and the L\'evy measure $\nu$ has density given by \eqref{e:stable-density}. Moreover, it holds that 
$d=a:=(c_+-c_-)/(\alpha-1)$ 
when $\alpha\neq 1$, and we specify 
$d=a=0$ when $\alpha=1$, cf.~\cite[p.~398]{KPW14}.

Recall also that $\tau=\tau_{(0,\infty)}:=\inf\{t>0: \eta_t\in (-\infty, 0]\}$, $X^{\ast}_t:=\eta_t\1_{(t<\tau)}$, $t\ge 0$, and
$\xi^{\ast}$ is the L\'evy process associated to $X^{\ast}$ through the Lamperti transform. Then $\xi^{\ast}$ is a killed L\'evy process. Its L\'evy measure  $\mu$  was computed in \cite[Section 3.1]{CC06}, see also \cite[(6)]{KPW14}. It holds that  $\mu$   has a
density
\begin{equation}\label{e:levy-xi*}
\mu(x) =c_+\frac{e^x}{(e^x-1)^{1+\alpha}}\1_{(x>0)}+
c_-\frac{e^x}{(1-e^x)^{1+\alpha}}\1_{(x<0)}, 
\quad x\in \R,
\end{equation} 
and $\xi^{\ast}$ is killed at rate
\begin{equation}\label{e:killing-rate}
\frac{c_-}{\alpha}=\frac{\Gamma(\alpha)}{\Gamma(\alpha\wh{\rho})\Gamma(1- \alpha\wh{\rho})}.
\end{equation}
The characteristic exponent $\Psi^{\ast}$  of $\xi^{\ast}$ can be found in \cite[(13.46)]{Kyp14}:
\begin{equation}\label{e:char-exp-xi*}
\Psi^{\ast}(\theta)=\frac{\Gamma(\alpha-i\theta)}{\Gamma(\alpha\wh{\rho}-i\theta)}\, \frac{\Gamma(i\theta+1)}{\Gamma(i\theta+1-\alpha\wh{\rho})}, \quad \theta \in \R.
\end{equation}
The infinitesimal generator of $\xi^*$ is given by
\begin{equation}\label{e:inf-ge-xi-star}
\AA^* f(x)=-\frac{c_-}{\alpha}f(x)-bf'(x)+\int_{\R}\left(f(x+y)-f(x)-f'(x)y \1_{[-1,1]}(y)\right) \mu(dy),
\end{equation}
where $b\in \R$ is a linear term which we will identify shortly.
By \eqref{e:pssMp-inf-gen-2}, the infinitesimal generator of $X^*$ is equal to 
\begin{eqnarray*}
\LL^* f(x)&=&
-\frac{c_-}{\alpha}x^{-\alpha}f(x)
 -b x^{1-\alpha}f'(x)+x^{-\alpha}\int_{\R}\left(f(xe^y)-f(x)-xf'(x)y\1_{[-1,1]}(y)\right)\mu(y)dy\\
&=&
-\frac{c_-}{\alpha}x^{-\alpha}f(x)
-b x^{1-\alpha}f'(x)\\
& & +x^{-\alpha}
\int^\infty_0
\left(f(z)-f(x)-xf'(x)(\log z/x)\1_{[-1,1]}(\log(z/x))\right)\mu(\log(z/x))z^{-1}dz.
\end{eqnarray*}
We have that
\begin{eqnarray*}
\mu\left(\log \frac{z}{x}\right)&=&c_+\frac{\frac{z}{x}}{\left(\frac{z}{x}-1\right)^{1+\alpha}}\1_{(z>x)}+c_-\frac{\frac{z}{x}}{\left(1-\frac{z}{x}\right)^{1+\alpha}}\1_{(z<x)}\\
&=& x^{\alpha} z \left(c_+ |z-x|^{-1-\alpha}\1_{(z>x)}+c_- |z-x|^{-1-\alpha}\1_{(z<x)}\right)\\
&=& x^{\alpha} z\, \nu(z-x).
\end{eqnarray*}
Therefore, with $j(x,z):=\nu(z-x)$, 
\begin{eqnarray}\label{e:gen-killed}
\LL^* f(x)&=&
-\frac{c_-}{\alpha}x^{-\alpha}f(x)
-b x^{1-\alpha}f'(x)\nonumber \\
 & &+
 \int^\infty_0
 \left(f(z)-f(x)-xf'(x)(\log z/x)\1_{[-1,1]}(\log(z/x))\right)j(x,z)dz.
\end{eqnarray}
By comparing this expression with the form of the infinitesimal generator of $X^*$ given in \cite[Theorem 2]{CC06}, we see that 
\begin{equation}\label{e:linear-term}
b= - a- 
\int^\infty_0
\left((\log u)\1_{[-1,1]}(\log u)-(u-1)\1_{[-1,1]}(u-1)\right)\nu(u-1)\, du.
\end{equation}

Let $\Psi(\theta):=\Psi^{\ast}(\theta)-\Psi^{\ast}(0)=\Psi^{\ast}(\theta)-c_- /\alpha$. Then $\Psi$ is the characteristic exponent of 
an unkilled  L\'evy processes $\xi$.
More precisely, $\xi^{\ast}$ is equal in distribution to $\xi$ killed at an independent exponential time with parameter $c_-/\alpha$. 
The infinitesimal generator of $\xi$ is
\begin{equation}\label{e:inf-gen-xi}
\AA f(x)= -bf'(x)+\int_{\R}\left(f(x+y)-f(x)-f'(x)y \1_{[-1,1]}(y)\right) \mu(dy).
\end{equation}
It is immediate from \eqref{e:levy-xi*} that 
$\int_{|y|\ge 1}|y|\mu(dy)<\infty$, hence $\E|\xi_1|<\infty$. 
Let $X=(X_t, \P_x)$ be the pssMp  of index $\alpha$
corresponding to $\xi$ through the Lamperti transform. The effect of removing the killing term from the generator of $\xi^{\ast}$ is to remove the killing term from the generator of $X^{\ast}$. Hence, the infinitesimal generator 
$\LL$ 
of $X$ is given by the right-hand side of \eqref{e:gen-killed} with 
$-(c_-/\alpha)x^{-\alpha}f(x)$ removed:
\begin{equation}\label{e:generator-LL-X}
\LL f(x)=-b x^{1-\alpha}f'(x)+\int^\infty_0\left(f(y)-f(x)-xf'(x)(\log y/x)\1_{[-1,1]}(\log(y/x))\right)j(x,y)dy.
\end{equation}
Considered on $(0,\infty)$, $X^{\ast}$ is a stable process killed upon exiting $(0,\infty)$. 
By removing 
the killing term in the infinitesimal generator, we end up with the process $X$ -- 
the  (not necessarily symmetric) \emph{censored $\alpha$-stable process} on $(0,\infty)$. 
The censored process $X$ can be also regarded as a resurrected
process with the resurrection kernel $q(x,A)=\P_x(\eta_{\tau_-}\in A)$ -- it is continued exactly at the position from which $\eta$ has jumped out from $(0,\infty)$ (thus effectively suppressing this jump). 
Note that this type of resurrection does not fall into our setting.
Censored processes were introduced in \cite{BBC} in a more general multi-dimensional context for rotationally symmetric stable processes.

\subsection{Examples of return kernels}\label{ss:examples-return}
An example of pssMp  of index $\alpha$ related to $\eta$ is its absolute value process $|\eta |=(|\eta_t|)_{t\ge 0}$. One can view 
$|\eta|$ also as a resurrected process: at time $\tau$, if $z=\eta_{ \tau}$, 
we resurrect at $-z>0$ according to the resurrection kernel $q$ with $p(z,A)=\delta_{-z}(A)$. 

We discuss now two families of return kernels with $\mathfrak{p}=1$.
First note that if the measure $\phi$ is absolutely continuous
 with respect to the Lebesgue measure on $[0, \infty)$ with a density (which we denote by the same letter), 
then the return kernel $p(z, \cdot)$ has a density $p(z,y)$, $y>0$, and \eqref{e:p-scaling-measure}--\eqref{e:p-phi-measure} imply that
$$
p(z,y)=\phi\left(\frac{y}{|z|}\right)\frac{1}{|z|}.
$$

In the first family of return kernels the probability measure $\phi$ has a density which decays polynomially at infinity:
For $\beta>0$ and $\gamma>\beta$, let
\begin{equation}\label{e:phi-special}
\phi(t)=\phi_{\beta,\gamma}(t)=
\frac{\Gamma(\gamma)}{\Gamma(\beta)\Gamma(\gamma-\beta)}t^{\beta-1}(1+t)^{-\gamma}.
\end{equation}
Motivation for this family comes from two particular examples. The first one is the path-censored process introduced in \cite{KPW14}. This process is obtained from 
$\eta$ by removing parts of the path in $(-\infty,0]$. More formally, define
$A_t:=\int_0^t \1_{(\eta_s>0)}ds$
and let $\tau_t:=\inf\{s>0:\, A_s>t\}$ be its right-continuous inverse. The process $\theta=(\theta_t)_{t\ge 0}$, defined by $\theta_t=\eta_{\tau_t}$,  is a strong Markov process 
on $(0,\infty)$, called the \emph{path-censored} process of $\eta$ on $(0,\infty)$.
The part of the process $\theta$ until its first hitting time of $0$ can be described in the following way: 
Let $x=\eta_{\tau-}\in (0,\infty)$ be the position from which $\eta$ jumps 
out of $(0,\infty)$, and $z=\eta_{\tau}<0$ be the position where $\eta$ lands at the exit from $(0,\infty)$. 
The distribution
of the returning position of $\eta$ to $(0,\infty)$ has the density $P_{(-\infty, 0)}(z,y)$ called the Poisson kernel: If $\sigma:=\inf\{t>0:\, \eta_t\in [0,\infty)\}$, then 
$\P_z(\eta_{\sigma}\in A)=\int_A P_{(-\infty,0)}(z,y)\, dy$, $A\in \BB((0, \infty))$.
The exact formula for this Poisson kernel is given by 
(e.g. \cite[Lemma 1.1]{Kyp18} which contains a minor typo: the $\alpha$ there should be $\alpha\rho$),
$$
P_{(-\infty,0)}(z,y)=\frac{1}{\Gamma(1-\alpha\rho)\Gamma(\alpha\rho)}\left(\frac{y}{|z|}\right)^{-\alpha\rho}(|z|+y)^{-1}=\phi_{1-\alpha\rho, 1}\left(\frac{y}{|z|}\right)\frac{1}{|z|}.
$$
In the second example the return kernel is equal to the normalized jump kernel, see \cite{DROV17, Von21}, 
\begin{align}
\label{e:njk}
p(z,y)=\frac{j(z,y)}{\int_0^{\infty}j(z,u)du}=\alpha |z|^{\alpha}(|z|+y)^{-1-\alpha}=\phi_{1,1+\alpha}\left(\frac{y}{|z|}\right)\frac{1}{|z|}.
\end{align}
Here we used that $j(z,y)=c_+(y-z)^{-1-\alpha}$ for $z<0$, $y>0$.

In the second family, the probability measure $\phi$ has a density which decays exponentially at infinity: 
For $a,\beta, \gamma>0$, let 
\begin{equation}\label{e:phi-special-exp}
\phi(t)=\phi_{a,\beta,\gamma}(t)=\frac{\gamma a^{\frac{\beta}{\gamma}}}{\Gamma(\frac{\beta}{\gamma})}t^{\beta-1} e^{-at^{\gamma}}, \quad t>0.
\end{equation}
Then
$$
p(z,y)=\frac{\gamma a^{\frac{\beta}{\gamma}}}{\Gamma(\frac{\beta}{\gamma})}\frac{y^{\beta-1}}{|z|^{\beta}} e^{-a\left({y}/{|z|}\right)^{\gamma}}.
$$


\section{Resurrected process}\label{s:res-proc}

Let us go back to the strictly $\alpha$-stable process $\eta$. 
Instead of killing $\eta$ upon exiting $(0,\infty)$ to get $X^{\ast}$, 
or restarting at $\eta_{\tau-}$ 
to get the censored process $X$, we look at the exit point 
$\eta_{\tau}=z<0$ 
and (partially)  resurrect
according to some probability kernel $p(z,A)$,   
$A\in \BB( [  0, \infty))$.
We assume that  $p(z, \cdot)$ satisfies \eqref{e:p-scaling-measure}. Recall from \eqref{e:p-phi-measure}  that such  kernels are in one-to-one correspondence with probability measures $\phi$ on  
$\BB([  0, \infty))$
through the relation
$$
p(z,A)=\phi(|z|^{-1}{A}), \quad z<0.
$$
If the measure $p(z,\cdot)$ has a density with respect to the Lebesgue measure, we will denote it by $p(z,y)$, 
$y>0$.
It is immediate that 
\begin{equation}\label{e:p-scaling}
p(\lambda z, \lambda y)=\lambda^{-1} p(z,y)\, \quad \text{for all }z<0, \text{ a.e. }y>0 \text{ and all }\lambda >0.
\end{equation}
Recall that $\nu(x)$ is the  L\'evy density given in \eqref{e:stable-density}. 
For $x,z\in \R$, let $j(x,z):=\nu(z-x)$.
If $x>0$ and $z<0$, we have that
\begin{equation}\label{e:j-nonsym}
j(x,z)=c_-(x-z)^{-1-\alpha}.
\end{equation}
Further, $j$ enjoys the following scaling property: 
\begin{align}
\label{e:spj}
j(\lambda x, \lambda z)=\lambda^{-1-\alpha}j(x,z), \quad x>0, z<0, \lambda >0.
\end{align}
Recall that the  resurrection kernel $q(x, A)$ was defined in \eqref{e:int-kernel-measure} as
\begin{align}
\label{e:qpph}
q(x,A):= \int_{ (-\infty, 0) }j(x,z)p(z,A)\, dz
=\int_{ (-\infty, 0)}\frac{c_-}{(x-z)^{1+\alpha}}\phi\left(\frac{A}{|z|}\right) dz
\quad x>0, A\in 
\BB((0, \infty)).
\end{align}

\subsection{Compound Poisson processes corresponding to resurrection  kernels
}\label{ss:cpp}

In this subsection we study the pssMp defined through the  resurrection kernel $q$ and the corresponding L\'evy process.
From the scaling properties of $p$ in 
\eqref{e:p-scaling-measure} 
 and $j$ in \eqref{e:spj} we get that the  
resurrection kernel $q$  satisfies
$$
q(\lambda x, \lambda A)=\lambda^{-\alpha}q(x,A), \quad x>0,  
A\in \BB((0, \infty)),\lambda >0.
 $$
In particular,
$q (1,A)=x^{\alpha }q(x,xA)$,
implying that for any $g:(0, \infty) \to \R$,
\begin{equation}\label{e:integral-1-x}
\int g(y)q(1,dy)=
x^{ \alpha}\int g(y/x)q(x,dy).
\end{equation}

\begin{lemma}\label{l:q-always-density}
For all $x>0$, 
the measure $q(x, \cdot)$ has a density $q(x,y)$ given by
\begin{equation}\label{e:q-always-density}
q(x,y):=c_-\int_{(0,\infty)} 
\left(x+\frac{y}{t}\right)^{-1-\alpha}t^{-1}\phi(dt).
\end{equation}
\end{lemma}
\pf 
By \eqref{e:qpph}, for any $A\in \BB((0, \infty))$ we have
\begin{eqnarray*}
q(x,A)&=&c_-
\int_0^\infty (x+z)^{-1-\alpha}\left(\int_{(0, \infty)}
\1_{A/z}(t)\phi(dt)\right) dz\\
&=&c_-
\int_{(0, \infty)} \int^\infty_0(x+z)^{-1-\alpha}\1_A(zt)\, dz\,  \phi(dt)\\
&=& c_-
\int_{(0, \infty)} \int^\infty_0\left(x+\frac{y}{t}\right)^{-1-\alpha}t^{-1}
\1_A(y)\, dy\, \phi(dt)\\
&=&\int_A \left(c_-
\int_{(0, \infty)}\left(x+\frac{y}{t}\right)^{-1-\alpha}t^{-1}\, \phi(dt)
\right) dy,
\end{eqnarray*}
where in the second and last equalities we used 
Tonelli's theorem  and in the penultimate equality  the change of variables $y=tz$.
\qed

We record here a simple consequence of \eqref{e:q-always-density}: For all $y>0$,
\begin{equation}\label{e:simple-bound-on-q}
q(1,y)\le c_- (y^{-1-\alpha}+1). 
\end{equation}
Indeed,
\begin{eqnarray*}
q(1,y)&=&c_-\left(\int_{(0,1)}\left(1+\frac{y}{t}\right)^{-1-\alpha}\frac{\phi(dt)}{t}+\int_{[1,\infty)}\left(1+\frac{y}{t}\right)^{-1-\alpha}\frac{\phi(dt)}{t}\right)\\
&\le & c_- \left( y^{-1-\alpha}\int_{(0,1)}t^{\alpha}\phi(dt)+\int_{[1,\infty)}t^{-1}\phi(dt)\right)\le c_- (y^{-1-\alpha}+1). 
\end{eqnarray*}

Note that
$$
q(x):=
q(x,(0, \infty))=\int_{-\infty}^0j(x,z)p(z,(0, \infty))dz
=\mathfrak{p}  c_-\int_{-\infty}^0(x-z)^{-1-\alpha}dz=  \mathfrak{p}  \frac{c_-}{\alpha}x^{-\alpha},
$$
so that
$$
Q(x,A):=\frac{q(x,A)}{q(x)}, \quad x>0, 
A\in \BB((0, \infty)),
$$
is a well-defined probability kernel satisfying
$$
Q(\lambda x, \lambda A)=Q(x,A), \quad x>0, 
A\in \BB((0, \infty)),
\lambda >0.
$$
It follows from Lemma \ref{l:q-always-density} that  both $q(x, \cdot)$ and $Q(x,\cdot)$  have densities $q(x,y)$, resp.~$Q(x,y)$,  satisfying
$$
q(\lambda x, \lambda y)=\lambda^{-1-\alpha}q(x,y), \quad Q(\lambda x, \lambda y)=\lambda^{-1}Q(x,y), \quad x,y>0, \lambda >0.
$$

We define now $\Pi(x,\cdot)$ to be the image measure of $Q(e^x,\cdot)$ under the mapping $y\mapsto e^y$: $\Pi(x,A):=Q(e^x, e^A)$. In particular, for every $g:\R\to \R$,
\begin{equation}\label{e:pi-q-change}
\int_{\R}g(y)\Pi(x,dy)=
\int_{(0,{\infty})}g(\log y)Q(e^x, dy).
\end{equation}
Note that $\Pi(\cdot, \cdot)$ is translation invariant, that is, for all $u\in \R$,
$$
\Pi(x+u, A+u)=Q(e^x e^u, e^u e^A)=Q(e^x, e^A)=\Pi(x,A).
$$
Let
$$
\Pi(A):=
\Pi(0,A)
=Q(1,e^A), \quad A\in \BB(\R),
$$
and 
$$
\pi(A):=q(1,e^A)=q(1)Q(1,e^A)=
q(1)\Pi(A), \quad A\in \BB(\R).
$$
Clearly, $\Pi$ is a probability measure on $\BB(\R)$ and $\pi$ a finite measure, hence a L\'evy measure. 
Both $\Pi$ and $\pi$ have densities (which again by an abuse of notation we denote by the same letters) satisfying
\begin{equation}\label{e:densities-Q-q}
\Pi(y)=Q(1,e^y)e^y \quad \text{and } \quad \pi(y)=q(1,e^y)e^y=e^{-\alpha y}q(e^{-y},1). 
\end{equation}
It follows from \eqref{e:simple-bound-on-q} that 
\begin{equation}\label{e:simple-bound-on-pi}
\pi(y)\le c_- (e^{-(1+\alpha)y}+1)e^y=c_- (e^{-\alpha y}+e^y).
\end{equation}
In particular, $\pi(y)$ is bounded in every neighborhood of 0.

Let $\chi$ be a compound Poisson process with characteristic exponent
$\Psi^\chi(\theta)=\int_{\R}(1-e^{i\theta y})\pi(dy)$. 
The jump distribution of $\chi$ is given by the measure $\Pi$ and 
the jump rate  is $q(1)$. 
The infinitesimal generator of $\chi$ is given by
$$
\AA f(x)=\int_{\R}(f(x+y)-f(x))\pi(dy).
$$

Let $Y=(Y_t, \P_x)$ be the pssMp of index $\alpha$
related to $\chi$ through  the Lamperti transform.
The infinitesimal generator of $Y$ is, according to \eqref{e:pssMp-inf-gen-2}, equal to
\begin{eqnarray}\label{e:new-number-for-LL}
\LL^Yf(x)&=&x^{-\alpha} \int_{\R}(f(xe^u)-f(x))\pi(du)=x^{-\alpha}
\int_{(0,{\infty})} (f(xu)-f(x))q(1,du) \nn \\
&=&
\int_{(0,{\infty})}(f(y)-f(x))q(x,dy),
\end{eqnarray}
where  we used \eqref{e:pi-q-change} in the second 
equality, and \eqref{e:integral-1-x} in the third equality.
This shows that $Y=(Y_t, \P_x)$ is a regular step process defined by the Markov kernel $Q(x,A)$, 
$A\in \BB((0, \infty))$, 
and the holding function $q(x)$, cf.~\cite[I.12]{BG68}.

\begin{thm}\label{t:ft-pi}
Suppose that $p(z,\cdot)$ is given by \eqref{e:p-phi-measure}. Then
$$
\hat{\pi}(\theta):=\int_{\R}e^{i\theta y}\pi(dy)=\frac{c_-}{\alpha}\frac{\Gamma(\alpha-i\theta)\Gamma(1+i\theta)}{\Gamma(\alpha)}
\int_{(0, \infty)}u^{i\theta}\phi(du).
$$
\end{thm}
\pf 
We have for $A\in \BB(\R)$,
\begin{eqnarray}
\pi(A)&=&q(1,e^A)=c_-
\int_0^\infty
(1+z)^{-1-\alpha}\phi\left(\frac{e^A}{z}\right)dz \nonumber \\
&=&c_- 
\int_0^\infty(1+z)^{-1-\alpha}\left(\int_{(0, \infty)}
\1_{e^A/z}(y)\phi(dy)\right)dz \nonumber \\
&=&c_- 
\int_0^\infty(1+z)^{-1-\alpha}\left(\int_{(0, \infty)}
\1_A(\log(yz))\phi(dy)\right)dz \nonumber \\
&=&c_- 
\int_{(0, \infty)}\left( \int_0^\infty
\1_A(\log(yz))(1+z)^{-1-\alpha}dz\right)\phi(dy). \label{e:pi}
\end{eqnarray}
Therefore, by using \cite[8.380.1-3]{GR07}, 
\begin{eqnarray*}
\hat{\pi}(\theta)
&=&\int_{\R}e^{i\theta y}\pi(dy)
= c_-
\int_{(0, \infty)}\left( \int_0^\infty
e^{i\theta\log(yz)}(1+z)^{-1-\alpha}dz \right)\phi(dy)\\
&=&c_- 
\int_{(0, \infty)}\left( \int_0^\infty
(yz)^{i\theta}(1+z)^{-1-\alpha}dz\right)\phi(dy)\\
&=& c_- \left(
\int_0^\infty z^{i\theta}(1+z)^{-1-\alpha}dz\right)
\left(\int_{(0, \infty)}y^{i\theta}\phi(dy)\right)\\
&=& \frac{c_-}{\alpha}\frac{\Gamma(\alpha-i\theta)\Gamma(1+i\theta)}{\Gamma(\alpha)}\, 
\int_{(0, \infty)}y^{i\theta}\phi(dy).
\end{eqnarray*}
\qed

This theorem allows
us to rewrite the density of the jump distribution of $\chi$  as a convolution of 
a subprobability and a probability distribution on $\R$,
cf.~\cite[p.411, 2nd paragraph]{KPW14}. Define 
$$
\tau(A):=\phi(e^A), \quad A\in \BB(\R), 
\quad\text{and} \quad f(y):=\frac{\alpha e^y}{(1+e^y)^{1+\alpha}}, \quad y\in \R.
$$
Then $f$ is a probability density on $\R$ and
$$
\hat{f}(\theta):=\int_{\R}e^{i\theta y}f(y)dy =\frac{\Gamma(\alpha-i\theta)\Gamma(1+i\theta)}{\Gamma(\alpha)}, \quad 
\hat{\tau}
(\theta):=\int_{\R}e^{i\theta y}\tau(dy)=
\int_{(0, \infty)}y^{i\theta}\phi(dy).
$$

\begin{corollary}\label{c:density-Pi}
It holds that 
$$
\Pi(y)=(f\ast \tau)(y)=\alpha e^{ y}
\int_{(0, \infty)}
\frac{t^\alpha}
{(t+e^{y})^{1+\alpha}}\phi(dt), \quad y\in \R.
$$
\end{corollary}
\pf The first equality is an immediate consequence of the equality $\hat{\Pi}=\hat{f}\hat{\tau}$. For the second equality, we rewrite
\begin{align*}\label{e:demsity-Pi-general}
\Pi(y)&= \int_{\R} \frac{\alpha e^{y-u}}{(1+e^{y-u})^{1+\alpha}}\tau(du)=\alpha e^y
\int_{(0, \infty)}
\frac{  t^{-1}}{(1+t^{-1}e^{y})^{1+\alpha}} \phi(dt) =\alpha e^{ y}
\int_{(0, \infty)}
\frac{t^\alpha}
{(t+e^{y})^{1+\alpha}}\phi(dt). \nonumber
\end{align*}
\qed

\begin{corollary}\label{c:ft-pi}
(a)  Assume that $\phi$ is given by \eqref{e:phi-special}. Then 
\begin{equation}\label{e:hat-pi}
\hat{\pi}(\theta)= \frac{c_-}{\alpha}\frac{\Gamma(\alpha-i\theta)\Gamma(1+i\theta)}{\Gamma(\alpha)} \frac{\Gamma(\beta+i\theta)\Gamma(\gamma-\beta-i\theta)}{\Gamma(\beta)\Gamma(\gamma-\beta)}  .
\end{equation}

\noindent 
(b) Assume that $\phi$ is given by \eqref{e:phi-special-exp}. Then
\begin{equation}\label{e:hat-pi-phi2}
\hat{\pi}(\theta)= \frac{c_-}{\alpha}\frac{\Gamma(\alpha-i\theta)\Gamma(1+i\theta)}{\Gamma(\alpha)}
\frac{a^{-\frac{i\theta}{\gamma}}\Gamma\left(\frac{\beta+i\theta}{\gamma}\right)
 }{\Gamma\left(\frac{\beta}{\gamma}\right)}.
\end{equation}

\noindent
(c) Assume that $\phi=\delta_a$, $a\in
(0,{\infty})$. Then
$$
\hat{\pi}(\theta)= \frac{c_-}{\alpha}\frac{\Gamma(\alpha-i\theta)\Gamma(1+i\theta)}{\Gamma(\alpha)} \, a^{i\theta}.
$$
\end{corollary}
\pf
(a) It follows from Theorem \ref{t:ft-pi} and \cite[8.380.1-3]{GR07}   that 
\begin{eqnarray*}
\hat{\pi}
(\theta)&=&\frac{c_-}{\alpha}\frac{\Gamma(\alpha-i\theta)\Gamma(1+i\theta)}{\Gamma(\alpha)} \frac{\Gamma(\gamma)}{\Gamma(\beta)\Gamma(\gamma-\beta)}\int_0^{\infty} \frac{u^{i\theta +\beta-1}}{(1+u)^{\gamma}}\, du\\
&=&\frac{c_-}{\alpha}\frac{\Gamma(\alpha-i\theta)\Gamma(1+i \theta)}{\Gamma(\alpha)} \frac{\Gamma(\beta+i\theta)\Gamma(\gamma-\beta-i\theta)}{\Gamma(\beta)\Gamma(\gamma-\beta)}  .
\end{eqnarray*}

\noindent (b) This follows from Theorem \ref{t:ft-pi} by 
noting 
(after the change of variables $v=au^{\gamma}$) that 
$$
\int_0^{\infty}u^{i\theta+\beta-1}e^{-au^{\gamma}} \, du= \gamma^{-1}a^{-\frac{\beta+i\theta}{\gamma}}\Gamma\left(\frac{\beta+i\theta}{\gamma}\right).
$$

\noindent (c) This is clear.
\qed

In case $\phi$ is given by \eqref{e:phi-special}, we can also compute  the density $\Pi(y)$  of the jump distribution of $\chi$. 
By using 
Corollary \ref{c:density-Pi} and the 
change of variables $u=t^{-1}e^y$
 we have 
\begin{eqnarray}\label{e:tilde-pi-exact}
\lefteqn{\Pi(y)
= 
\frac{\alpha\Gamma(\gamma)}{\Gamma(\beta)\Gamma(\gamma-\beta)} e^y \int_0^{\infty} \left(t+e^y\right)^{-1-\alpha}t^{\beta+\alpha-1}(1+t)^{-\gamma}dt }\nonumber \\
&=& \frac{\alpha\Gamma(\gamma)}{\Gamma(\beta)\Gamma(\gamma-\beta)} e^{-(\gamma-\beta) y}\int_0^{\infty} u^{\gamma-\beta}(1+ue^{-y})^{-\gamma}(1+u)^{-1-\alpha}du \nonumber \\
&=&\frac{\alpha\Gamma(\gamma)B(1+\gamma-\beta, \beta-\alpha)}{\Gamma(\beta)\Gamma(\gamma-\beta)} e^{-(\gamma-\beta) y}\,  {_2}\FF_1(\gamma, 1+\gamma-\beta; 1+\alpha+\gamma; 1-e^{-y}),
\end{eqnarray}
where the last line follows from \cite[3.197.5]{GR07}. Here ${_2}\FF_1$ is the hypergeometric function.

\begin{example}\label{ex:trace}{\rm
(a) Recall that $\alpha \rho<1$. Take $\beta=1-\alpha \rho$ and $\gamma=1$ in \eqref{e:phi-special}, we get
\begin{equation}\label{e:res-trace}
p(z,y)=\frac{1}{\Gamma(1-\alpha \rho)\Gamma(\alpha \rho)}\frac{|z|^{\alpha \rho}}{y^{\alpha\rho}}(|z|+y)^{-1}.
\end{equation}
Since $\Pi(y)=(\alpha/c_-)\pi(y)$, it follows from Corollary \ref{c:ft-pi} that the characteristic function of the jump distribution $\Pi$ of $\chi$ is equal to
$$
\frac{\Gamma(\alpha-i\theta)\Gamma(1+i\theta)}{\Gamma(\alpha)} \frac{\Gamma(1-\alpha\rho+i\theta)\Gamma(\alpha \rho-i\theta))}{\Gamma(1-\alpha \rho)\Gamma(\alpha \rho)}  .
$$
Since $\Gamma(1-\alpha \rho)\Gamma(\alpha \rho)=\frac{\pi}{\sin(\pi\alpha\rho)}$,  cf.~\cite[8.334.3]{GR07}, 
we see that the characteristic function of the jump distribution of $\chi$ coincides with the one in  \cite[Proposition 4.2, (9)]{KPW14}. This means that the  return kernel given in \eqref{e:res-trace} corresponds to the path censored (or trace) process studied in \cite{KPW14}. This can be also seen by recognizing $p(z,y)$ from \eqref{e:res-trace} as the Poisson kernel $P_{(-\infty,0)}(z,y)$ of the stable process 
$\eta$
(see Subsection \ref{ss:examples-return} above). 
By using \eqref{e:tilde-pi-exact}  with $\beta=1-\alpha\rho$ and $\gamma=1$, 
we find
$$
\Pi(y)=\frac{\alpha\Gamma(\alpha \rho+1)\Gamma(\alpha\wh{\rho}+1)}{\Gamma(1-\alpha\rho)\Gamma(\alpha\rho)\Gamma(\alpha+2)}
 \, e^{-\alpha\rho y}\,  _2\FF_1(1, \alpha\rho+1; \alpha+2; 1-e^{-y}),\quad y\in \R,
$$
cf.~\cite[(13)]{KPW14}.

\noindent 
(b) Take $\beta=1$ and $\gamma=1+\alpha$ in \eqref{e:phi-special}. Then
\begin{equation}\label{e:rez-Z}
p(z,y)=\frac{\Gamma(1+\alpha)}{\Gamma(\alpha)}|z|^{\alpha}(|z|+y)^{-1-\alpha}=\alpha|z|^{\alpha}(|z|+y)^{-1-\alpha},
\end{equation}
which is the normalized jump kernel in \eqref{e:njk}.
It follows from Corollary \ref{c:ft-pi} that the characteristic function of the jump distribution of $\chi$ is equal to
$$
\hat{\Pi}
(\theta)=\frac{\Gamma(\alpha-i\theta)^2 \Gamma(1+i\theta)^2}{\Gamma(\alpha)^2},
$$
and the density is
$$
\Pi(y)=\alpha^2 B(1+\alpha, 1-\alpha)e^{-\alpha y} \,  _2\FF_1(1+\alpha, 1+\alpha; 2+2\alpha; 1-e^{-y}),\quad y\in \R.
$$

\noindent (c) Suppose that $\phi=\delta_a$, $a>0$.  Then
$$
\Pi(y)=
\alpha e^{ y}
\frac{a^\alpha}
{(a+e^{y})^{1+\alpha}}, \quad y\in \R.
$$
}
\end{example}

We end this subsection with a necessary and sufficient condition for  $\chi_1$  to have finite expectation.

\begin{prop}\label{p:chi-exp}
${\bf E}|\chi_1|<\infty$  if and only if 
\eqref{e:phi-int-log} holds true. 
\end{prop}
\pf 
It follows from \cite[Theorem 25.3, Example 25.12]{Sat14} that
${\bf E}|\chi_1|<\infty$ if and only if
$\int_{|y|\ge 1}|y|\pi(dy)<\infty$. 

We first assume \eqref{e:phi-int-log}.
We will show that $\int_{\R}|y|\pi(dy)<\infty$. By using \eqref{e:pi} we have
\begin{eqnarray*}
\int_{\R}|y|\pi(dy)&=&c_-
\int_{(0, \infty)}\left(\int_0^\infty
|\log(yz)|(1+z)^{-1-\alpha}dz \right)\phi(dy) \\
&\le & c_-
\int_{(0, \infty)}\left(\int_0^\infty
(|\log y| +|\log z|)|(1+z)^{-1-\alpha}dz \right)\phi(dy) \\
&=& c_-
\int_{(0, \infty)}|\log y| \left(\int_0^\infty
(1+z)^{-1-\alpha}dz\right)\phi(dy)\\
 & & + c_-
 \int_{(0, \infty)}\left(\int_0^\infty
 |\log z|(1+z)^{-1-\alpha}dz\right)\phi(dy) <\infty.
\end{eqnarray*}
Finiteness follows from the assumption and the fact that both integrals with respect to $dz$ are finite.

We now assume that 
$\int_{|y|\ge 1}|y|\pi(dy)<\infty$.
Then by \eqref{e:pi} 
\begin{align*}
\infty&>\int_{y \le -1}(-y)\pi(dy) \ge c_-
\int_{(0, 1)}\left(\int_0^\infty 
(-\log(yz)) (1+z)^{-1-\alpha} {\bf 1}_{\log(yz)\le -1} 
dz \right)\phi(dy) \\
& \ge c_-
\int_{(0, 1)}\left(\int_0^\infty 
 (-\log y)(1+z)^{-1-\alpha}{\bf 1}_{\log z\le -1} 
dz \right)\phi(dy)\\
&= c_-\left(\int_0^{1/e}
(1+z)^{-1-\alpha}dz\right)
\int_{(0, 1)}(-\log y) \phi(dy).
\end{align*}
Thus, $\int_{(0, 1)}(-\log y) \phi(dy)< \infty.$
Similarly, 
\begin{align*}
\infty&>\int_{y \ge 1}y\pi(dy) \ge c_-
\int_{(1, \infty)}\left(\int_0^\infty 
\log(yz)(1+z)^{-1-\alpha} {\bf 1}_{\log(yz)\ge 1}  \, dz \right)\phi(dy) \\
& \ge c_- \int_{(1, \infty)}\left(\int_0^\infty 
\log y(1+z)^{-1-\alpha} {\bf 1}_{\log z\ge 1} \, dz \right)  \phi(dy) \\
& = c_-\left(\int_{e}^\infty
(1+z)^{-1-\alpha}dz\right)
\int_{(1, \infty)}(\log y) \phi(dy).
\end{align*}
Thus, $\int_{(1, \infty)} ( \log y  ) \phi(dy)<\infty$. We have shown that \eqref{e:phi-int-log} holds. 
 \qed

\subsection{Resurrected process}\label{ss:res-proc}
Let $\Psi(\theta)=\Psi^{\ast}(\theta)-\Psi^{\ast}(0)=\Psi^{\ast}(\theta)-c_-/\alpha$ be the characteristic exponent of the L\'evy process $\xi$ corresponding to the censored $\alpha$-stable process $X$ through  the Lamperti transform. 
We define
$$
\Psi^{\sharp}(\theta):=\Psi(\theta)+(1-\mathfrak{p})\frac{c_-}{\alpha}=\Psi^{\ast}(\theta)-\mathfrak{p}\frac{c_-}{\alpha},
$$
and let $\xi^{\sharp}$ be the (killed) L\'evy process with the characteristic exponent $\Psi^{\sharp}$. 
 We will add to  $\xi^{\sharp}$   an independent compound Poisson process, denoted by $\chi$, with  characteristic exponent $\Psi^{\chi}$
given by
\begin{equation}\label{e:ch-exp-CPP}
\Psi^{\chi}(\theta)=\int_{\R}(1-e^{i\theta y})\pi(dy) = \mathfrak{p} \frac{c_-}{\alpha}-\hat{\pi}(\theta), \quad \theta\in \R.
\end{equation}
The effect of this procedure is that, 
instead of completely removing the killing from $\xi^{\ast}$,
we remove part of the killing (i.e., $\mathfrak{p}c_-/\alpha$) and add jumps according to $\Pi$ at the exponential rate 
$\wh{\pi}(0)=\mathfrak{p}c_-/\alpha$. 
Let $\overline{\xi}:= \xi^{\sharp} +\chi$ be this new L\'evy process. 
By \eqref{e:char-exp-xi*}, its characteristic exponent $\overline{\Psi}$ is given by
\begin{align}\label{e:char-exp-final}
\overline{\Psi}(\theta) = \Psi^{\sharp}(\theta)+\Psi^{\chi}(\theta) 
=\frac{\Gamma(\alpha-i\theta)}{\Gamma(\alpha\wh{\rho}-i\theta)}\, \frac{\Gamma(i\theta+1)}{\Gamma(i\theta+1-\alpha\wh{\rho})}-
\hat{\pi}
(\theta).
\end{align}
Note that 
$$
\overline{\Psi}(0)=\frac{\Gamma(\alpha)}{\Gamma(\alpha\wh{\rho})\Gamma(1-\alpha\wh{\rho})}-\wh{\pi}(0)=\frac{c_-}{\alpha}(1-\mathfrak{p}).
$$
If $\mathfrak{p}=1$, then $\overline{\Psi}(0)=0$, implying that in the Lamperti trichotomy the case (3) does not occur. If $\mathfrak{p}<1$, $\overline{\xi}$ is a killed L\'evy process 
with rate $\overline{\Psi}(0)$ and thus ${\bf E}{\overline \xi}_1=-\infty$.

Let $\overline{X}=(\overline{X}_t, \P_x)$ be the pssMp of index $\alpha$ with origin as a trap 
corresponding to $\overline{\xi}$ through the Lamperti transform. This process can be described as follows: Consider the strictly $\alpha$-stable process $\eta=(\eta_t, \P_x)$, and set as before $\tau=\inf\{t>0: \eta_t\le 0\}$. Then $\tau<\infty$ a.s. At time $\tau$,  
with probability $\mathfrak{p}$ we resurrect according to the  return kernel $p(\eta_{\tau}, y)$, $y>0$, and with probability $1-\mathfrak{p}$ kill the process and send it to the origin. 
This amounts to adding the  resurrection kernel $q(x,y)$ to the jump kernel $j(x,y)$ of $\eta$. 
The process $\overline{X}$ is a pssMp of index $\alpha$ with origin as a trap 
and jump kernel $J(x,y):=j(x,y)+q(x,y)$. 
Indeed, let $\AA^{\xi}$ be the infinitesimal generator of $\xi$ and $\LL^X$ the infinitesimal generator of $X$
given by \eqref{e:generator-LL-X}. 
The infinitesimal generator of $\overline{\xi}$ is obtained by adding $\pi$ to the L\'evy measure 
$\nu$ of $\xi$, and taking into account the killing term. 
Hence for $g\in C_0^2(\R)$, 
$$
\AA^{\overline{\xi}}g=-\frac{c_-}{\alpha}(1-\mathfrak{p})g+\AA^{\xi}g+\AA^{\chi}g.
$$
By using the relation between generators of the L\'evy process and the corresponding the Lamperti transformed pssMp of index $\alpha$, 
together with \eqref{e:new-number-for-LL},
we see that the infinitesimal generator of $\overline{X}$ is equal to
$$
\LL^{\overline{X}}f(x)=  - \frac{c_-}{\alpha}(1-\mathfrak{p})x^{-\alpha}+  \LL^X f(x)+\int_0^{\infty}(f(y)-f(x))q(x,y)dy.
$$

\medskip
Assume that $p(z,y)$ is given by \eqref{e:p-phi-measure}. 
Then by 
\eqref{e:c+c-}, \eqref{e:char-exp-final} and Theorem \ref{t:ft-pi},
\begin{align}\label{e:overline-Psi-general}
&\overline{\Psi}(\theta)=\frac{\Gamma(\alpha-i\theta)\Gamma(i\theta+1)}{\Gamma(\alpha\wh{\rho}-i\theta)\Gamma(i\theta+1-\alpha\wh{\rho})}\, 
-\frac{c_-}{\alpha}\frac{\Gamma(\alpha-i\theta)
\Gamma(i\theta+1)}{\Gamma(\alpha)}
\int_{(0,{\infty})}u^{i\theta}\phi(du)\nonumber \\
&=\Gamma(\alpha-i\theta)\Gamma(i\theta+1)\left(\frac{1}{\Gamma(\alpha\wh{\rho}-i\theta)\Gamma(i\theta+1-\alpha\wh{\rho})}-\frac{1}{\Gamma(\alpha\wh{\rho})\Gamma(1-\alpha\wh{\rho})}
\int_{(0,{\infty})}u^{i\theta}\phi(du)\right)\nonumber\\
&=\frac{\Gamma(\alpha-i\theta)
\Gamma(i\theta+1)}{\pi} \left(\sin(\pi(\alpha\wh{\rho}-i\theta)) -\sin(\pi\alpha\wh{\rho})
\int_{(0,{\infty})}u^{i\theta}\phi(du)\right).
\end{align}
In the third line we used the 
identity $\Gamma(z)\Gamma(1-z)=\pi/\sin (\pi z)$ twice.
\begin{remark}\label{r:ricochet}
{\rm
 For the ricocheted stable process from \cite{KPV21}, the measure $\phi$ determining the return kernel is equal to 
$(1-\mathfrak{p})\delta_0+\mathfrak{p}\delta_1$. 
In this case
$$
\overline{\Psi}(\theta)=\frac{\Gamma(\alpha-i\theta)
\Gamma(i\theta+1)}{\pi} \left(\sin(\pi(\alpha\wh{\rho}-i\theta)) -\mathfrak{p}\sin(\pi\alpha\wh{\rho})\right)
$$
which recovers \cite[(4.2)]{KPV21}, 
}
\end{remark}

\subsection{Behavior of $\overline{X}$ at absorption time }\label{ss:X-lifetime} 
If $\mathfrak{p}<1$, it follows from the Lamperti trichotomy that case (3) occurs, hence $\overline{X}$ is absorbed 
at 0 by a jump.
In the remaining part of this subsection we therefore assume that $\mathfrak{p}=1$. 
Recall that ${\bf E}|\xi_1|<\infty$ and, under assumption \eqref{e:phi-int-log}, also ${\bf E}|\chi_1|<\infty$, 
cf.~Proposition \ref{p:chi-exp}. 
Therefore, under assumption \eqref{e:phi-int-log}, 
${\bf E}|\overline{\xi}_1|<\infty$, $\overline{\Psi}'(0)$ exists, and ${\bf E}[\overline{\xi}_1]=i\overline{\Psi}'(0)$. 
Thus, combining  \cite[Theorem 7.2]{Kyp14} with the Lamperti trichotomy, we get that
if $i\overline{\Psi}'(0) \ge 0$, then $\limsup_{t\to \infty}\overline{\xi}_t=+\infty$, hence the 
absorption time of $\overline{X}$ is infinite; and if 
$i\overline{\Psi}'(0)< 0$, then $\lim_{t\to \infty}\overline{\xi}_t=-\infty$, hence the 
absorption time of $\overline{X}$ is finite $\P_x$-a.s.~and  
$\overline{X}$ is continuously absorbed at 0.

\medskip
\noindent
\textbf{Proof of Theorem \ref{t:derivative-at-zero}:}
The equivalence of $\mathbf{E}|\overline{\xi}_1|<\infty$  and \eqref{e:phi-int-log} follows from Proposition \ref{p:chi-exp}.
Put
$$
f_1(\theta):=B(\alpha-i\theta, 1+i\theta),\quad f_2(\theta):=
\int_{(0,{\infty})}u^{i\theta}\phi(du)
$$
$$
f_3(\theta)
:=\frac{\sin(\pi(\alpha\widehat\rho-i\theta))}{\pi} -\frac{\sin(\pi \alpha\widehat\rho \, )}{\pi}f_2(\theta),
$$
where $B$ denotes the beta function. 
Then by \eqref{e:overline-Psi-general}, 
$
\overline\Psi(\theta)=\Gamma(1+\alpha)f_1(\theta)f_3(\theta).
$
Since $f_2(0)=1$ and $f_3(0)=0$, we have 
\begin{align}
\label{oP0'}
\overline\Psi'(0)=\Gamma(1+\alpha)(f'_1(0)f_3(0)+f_1(0)f'_3(0))=
\Gamma(1+\alpha)B(\alpha, 1)f'_3(0)=
\Gamma(\alpha)f'_3(0).
\end{align}
Using the assumption \eqref{e:phi-int-log},
$$
f'_2(0)=i \int_{(0,{\infty})}(\log u)u^{i\theta}\phi(du)|_{\theta=0}=i\int_{(0,{\infty})}(\log u)\phi(du),
$$
and so
$$
f'_3(0)=-i\cos(\pi\alpha\widehat\rho)-\frac{\sin(\pi \alpha\widehat\rho)}{\pi}f'_2(0)=
-i\cos(\pi\alpha\widehat\rho)-i\frac{\sin(\pi \alpha\widehat\rho)}{\pi}
\int_{(0,{\infty})}(\log u)\phi(du).
$$
Therefore  by \eqref{oP0'}
\begin{align}\label{e:i-Psi-prime}
i\overline\Psi'(0)
&=\Gamma(\alpha)\left(\cos(\pi\alpha\widehat\rho)+\frac{\sin(\pi \alpha\widehat\rho)}{\pi}
\int_{(0,{\infty})}(\log u)\phi(du)\right)\nonumber\\
&=\Gamma(\alpha)\frac{\sin(\pi \alpha\widehat\rho)}{\pi}
\left(\pi\cot(\pi\alpha\widehat\rho)+
\int_{(0,{\infty})}(\log u)\phi(du)\right). \nonumber 
\end{align}
\qed

Since $(\alpha-1)_+ < \alpha\wh{\rho}<1$, we see from the display above that the sign of $i\overline{\Psi}'(0)$ depends on the sign of 
$$
\pi\cot(\pi\alpha\widehat\rho)+
\int_{(0,{\infty})}(\log u)\phi(du).
$$
Note that $\int_{(0,{\infty})}(\log u)\phi(du)$ may depend on $\alpha$. 
For example, when $\phi$ is given by \eqref{e:phi-special} with $\beta=1$ and $\gamma=1+\alpha$, 
it holds that $\int_{(0,{\infty})}(\log u)\phi(du)=\psi(1)-\psi(\alpha)$, 
where $\psi$ is the digamma function. 

Let $L_\phi:=-
\int_{(0,{\infty})}(\log u)\phi(du)$. 
Define $\arccot: \R\to (0,\pi)$ as a strictly decreasing and continuous function.
Set  
$$
a_{\phi}:
=\frac{1}{ \pi}\arccot\left(
\frac{L_\phi}{\pi}\right)
$$
and note that 
$a_{\phi}\in (0,1)$. 
\begin{corollary}\label{c:s}
Suppose $\mathfrak{p}=1$. 
(a) If  
$\alpha \le 1+a_{\phi}$,
 then  $\overline{X}$ is (continuously) absorbed at $0$ at an a.s.-finite time if and only if $\alpha \wh{\rho}>a_{\phi}$.
(b) If 
$\alpha>1+a_{\phi}$, 
then the 
absorption time of $\overline{X}$ is always finite $\P_x$-a.s..
\end{corollary}
\pf
(a)
If $\alpha \le 1$, since $0< \alpha\wh{\rho}<1$, we see that 
$\pi\cot(\pi\alpha\widehat\rho)-L_\phi<0$ 
if and only if
$$
\widehat\rho>
\frac{1}{ \alpha\pi}\arccot 
\left(\frac{L_\phi}{\pi}\right)=\frac{a_{\phi}}{ \alpha}.
$$
If $1<\alpha \le 1+a_{\phi}$, then 
$$
\cot(\pi(\alpha-1)) \ge 
\cot(\pi a_{\phi})= 
\frac{L_\phi}{\pi}
$$
and so we also have that 
$\pi\cot(\pi\alpha\widehat\rho)-
L_\phi<0$ 
if and only if $\widehat\rho>\frac{a_{\phi}}{ \alpha}$.

\noindent
(b) 
If $\alpha > 1+a_{\phi}$ 
then $\cot(\pi(\alpha-1)))< \frac{L_\phi}{\pi}$ 
and we always have
$$
\pi\cot(\pi\alpha\widehat\rho)-
L_\phi < \pi\cot(\pi(\alpha-1))-
L_\phi \le 0.
$$
\qed

In case $L_\phi=-\int_{(0,{\infty})}(\log u)\phi(du)$ is independent of $\alpha$, we can be slightly more precise. Let
$$
\rho(\alpha):=1-\frac{1}{\alpha \pi}\arccot\left(\frac{L_\phi}{\pi}\right),
$$
so that  
$\pi \cot\big(\pi \alpha(1-\rho(\alpha))\big)-
L_\phi=0$.   
Notice that $\rho(a_{\phi})=0$ (this need not be true in case $L_\phi$ depends on $\alpha$). 
Since $L_\phi$ does not depend on $\alpha$, the function $\alpha\mapsto \rho(\alpha)$ is strictly increasing, hence $a_{\phi}$ is the only zero of $\rho(\alpha)$.

\begin{corollary}\label{c:s2}
Suppose that $\mathfrak{p}=1$ and $L_\phi$ does not depend on $\alpha$. 
\begin{itemize}
\item[(i)] If  $\alpha\in (0, a_{\phi})$, then  the absorption time of $\overline{X}$ is  infinite;
\item[(ii)] If $\alpha\in [a_{\phi}, 1+a_{\phi}]$, then   the absorption time of $\overline{X}$ is finite if and only if $\rho>\rho(\alpha)$;
\item[(iii)] If $\alpha\in (1+a_{\phi}, 2)$, then the absorption time of $\overline{X}$ is finite.
\end{itemize}
\end{corollary}
\pf This is a direct consequence of Corollary \ref{c:s} and the discussion above. 
\qed

\subsection{Recurrent extension}\label{ss:rec-ext}
Recall that the origin is a trap for $\overline{X}$.  
If $\overline{X}$ is absorbed in 0 at finite time, 
one can ask if there exists  a positive self-similar  recurrent extension 
of $\overline{X}$. The general result is given in \cite[Theorem 1]{Fit06} and \cite[Theorems 1 and 2]{Riv07}: 
(i) 
There exists a unique 
positive self-similar  recurrent extension 
of $\overline{X}$ which leaves 0 continuously if and only if there exists $\kappa \in(0,\alpha)$ such that 
\begin{equation}\label{e:cont-rec-ext}
{\bf E}\left[e^{\kappa \overline{\xi}_1}\right]=1,
\end{equation}
and (ii) For $\beta\in(0,\alpha)$ there exists a positive self-similar recurrent extension of $\overline{X}$ which leaves 0 by a 
jump associated with an excursion measure of 
the form $c\beta x^{-1-\beta} dx$ if and only if 
\begin{equation}\label{e:jump-rec-ext}
{\bf E}\left[e^{\beta \overline{\xi}_1}\right]<1.
\end{equation}

Note that 
${\bf E}[e^{\kappa \overline{\xi}_1}]={\bf E}[e^{i(-i\kappa)\overline{\xi}_1}]=e^{-\overline{\Psi}(-i\kappa)}$ 
for all $\kappa\ge 0$ for which the expectation is finite. Let $\varphi:\R\to (-\infty, +\infty]$ be defined by $\varphi(\kappa):=-\overline{\Psi}(-i\kappa)$, so that 
$$
{\bf E}\left[e^{\kappa \overline{\xi}_t}\right]=e^{t\varphi(\kappa)}.
$$
Hence, \eqref{e:cont-rec-ext} is equivalent to the existence of $\kappa\in (0,\alpha)$ such that $\varphi(\kappa)=0$,
and \eqref{e:jump-rec-ext} is equivalent to $\varphi(\beta)<0$. 
Note that,  by H\"older's inequality, $\varphi$ is convex.

\medskip

\noindent
\textbf{Proof of Theorem \ref{t:R_e}:}
Let
$$
h(\kappa):=\sin(\pi(\alpha\wh{\rho}-\kappa))-\sin(\pi\alpha\wh{\rho})
\int_{(0,{\infty})}u^{\kappa}\phi(du).
$$
Clearly, $h(0)=(1-\mathfrak{p})\sin(\pi\alpha\wh{\rho})\ge 0$ since $\alpha\wh{\rho}\in (0,1)$,  and 
note that  from \eqref{e:overline-Psi-general} 
\begin{align}\label{e:Psik}
-\varphi(\kappa)=\overline{\Psi}(-i\kappa)
&=\frac{\Gamma(\alpha-\kappa)\Gamma(\kappa+1)}{\pi}
h(\kappa).
\end{align}

If  $\kappa_0<\alpha$, 
then $\kappa_0+\epsilon <\alpha$ for  all small $\epsilon >0$.
Since $\int_{(0,\infty)}u^{\kappa_0+\epsilon}\phi(du)=+\infty$ by definition of $\kappa_0$, we get that $h(\kappa_0+\epsilon)=-\infty$.

Assume that \eqref{e:kappa_0} holds true, i.e., $\kappa_0>0$. 
If $\kappa_0 \ge\alpha$, then 
\begin{eqnarray*}
h(\alpha-)&=&\sin(\pi(\alpha\wh{\rho}-\alpha))-\sin(\pi\alpha\wh{\rho})\lim_{\kappa\uparrow \alpha}
\int_{(0,{\infty})}u^{\kappa}\phi(du) \\
&=&-\sin(\pi\alpha \rho)-\sin(\pi\alpha\wh{\rho})\lim_{\kappa\uparrow \alpha}
\int_{(0,{\infty})}u^{\kappa}\phi(du) <0,
\end{eqnarray*}
where in the last inequality we used the assumptions $\alpha\rho\in (0, 1)$ and $\alpha\wh{\rho}\in (0, 1)$.
Therefore, $h((\alpha\wedge (\kappa_0+\epsilon)) -)<0$ for all small $\epsilon >0$.
If $\kappa\in  (  0,\kappa_0)$,
\begin{align}
\label{e:hprime}
h'(\kappa)=-\pi\cos(\pi(\alpha\wh{\rho}-\kappa))-\sin(\pi\alpha\wh{\rho})\left(
\int_{(0,{\infty})}(\log u)u^{\kappa}\phi(du) \right), 
\end{align}
which is justified by \eqref{e:kappa_0}. 

Assume that $\mathfrak{p}=1$. Since by the assumption  that $\overline{X}$ is absorbed at 0 in finite time, this happens continuously, and therefore 
${\bf E}[\overline{\xi}_1] \in [-\infty, 0)$. 
If \eqref{e:phi-int-log} holds true, since $i\overline{\Psi}'(0)<0$,  by Theorem \ref{t:derivative-at-zero} and \eqref{e:hprime},
$$
h'(0+)=-\pi\cos(\pi\alpha\wh{\rho})-\sin(\pi\alpha\wh{\rho}) 
\int_{(0,{\infty})}(\log u)\phi(du)>0,
$$
 implying that $h$ is strictly positive in some neighborhood of zero.
 Note that by \eqref{e:kappa_0}, 
we have  $\int_{(1,{\infty})}(\log u)u^{\kappa}\phi(du)<\infty$ for $\kappa\in (  0,\kappa_0)$ and 
thus $\int_{(1,{\infty})}(\log u)\phi(du)\le \int_{(1,{\infty})}(\log u)u^{\kappa}\phi(du)<\infty$. Consequently, if \eqref{e:phi-int-log} does not holds, then
we have  $\int_{(0,1)}(\log u)\phi(du)=-\infty$. 
 By the monotone convergence theorem, 
$$ \lim_{\kappa \downarrow 0}
 \int_{(0,1)}(\log u)u^{\kappa}\phi(du)= 
 -\lim_{\kappa \downarrow 0}\int_{(0,1)}(-\log u)u^{\kappa}\phi(du)=-\infty.
$$
We now see from \eqref{e:hprime} that 
$$
\liminf_{\kappa \downarrow 0}
h'(\kappa)=-\pi\cos(\pi\alpha\wh{\rho})-\sin(\pi\alpha\wh{\rho}) \left(\lim_{\kappa \downarrow 0}
\int_{(0,{\infty})}(\log u)u^{\kappa}\phi(du) \right)
=\infty,
$$
implying again that $h$ is strictly positive in some neighborhood of zero.
 
If $\mathfrak{p}<1$, then $h(0)>0$, so again we see that $h$ 
s strictly positive in some neighborhood of zero. 

Together with  
$h((\alpha\wedge (\kappa_0+\epsilon)) -)<0$ for all small $\epsilon >0$,
this implies the existence of $\kappa^{\ast}\in (0, \alpha)$ such that $h(\kappa^{\ast})=0$, hence also $\varphi(\kappa^{\ast})=0$.
Thus $\overline{X}$ has a positive self-similar recurrent extension.
Furthermore, by the convexity of $\varphi$, for every $\beta\in (0,\kappa^{\ast})$ we have $\varphi(\beta)<0$.
This means that 
${\bf E}[e^{\beta \overline{\xi}_1}]<1$.
By \cite[Theorem 1]{Riv07}, there exists a positive self-similar recurrent extension of $\overline{X}$ which leaves 0 
by a jump associated with the excursion measure $c\beta x^{-1-\beta}dx, x>0$.

Assume that \eqref{e:kappa_0} is false, that is $\kappa_0 =  0$. Then $h(\epsilon)=-\infty$ for all $\epsilon>0$, 
and we see from \eqref{e:Psik} that $\varphi(\beta)=+\infty$ for all $\beta\in (0,\alpha)$ and consequently 
${\bf E}[e^{\beta \overline{\xi}_1}]=+\infty$. 
Hence $\overline{X}$ does not have a recurrent extension.
\qed 
 
\begin{remark}\label{r:rec-ext-jump}
{\rm
It is easy to find examples of probability measures $\phi$ on $(0,\infty)$ satisfying \eqref{e:phi-int-log} but not \eqref{e:kappa_0} giving rise to pssMp that are 
(continuously) 
absorbed at zero in finite time, but not having a positive self-similar recurrent extension. One such example is the measure with density $\1_{(2,\infty)}(t)\frac{2(\log 2)^2}{t(\log t)^3}$.
}
\end{remark}

\subsection{Examples}\label{ss:examples}
In this subsection we analyze a list of examples.
Recall that  $\psi$ is the digamma function. 

\begin{example}\label{ex:beta-gamma}
{\rm
We look at our main example in which  $\phi$ is given by \eqref{e:phi-special}. 
 Let $f(t)=(1+t)^{-\gamma}$ with $\gamma>0$.
 The Mellin transform of $f$ is, by \cite[p.1131, 17.43.7]{GR07}, equal to
$$
M_f(s):=\int_0^{\infty}f(t)t^{s-1}dt=\frac{\Gamma(s)\Gamma(\gamma-s)}{\Gamma(\gamma)}.
$$
Since 
$$
\int_0^{\infty}f(t)(\log t)t^{s-1} dt=M_f'(s)=
\frac{1}{\Gamma(\gamma)}
\Big(\Gamma'(s)\Gamma(\gamma-s)-\Gamma(s)\Gamma'(\gamma-s)\Big),
$$
we get that
$$
\int_0^{\infty}(\log t)\phi(t)dt =\frac{\Gamma(\gamma)}{\Gamma(\beta)\Gamma(\gamma-\beta)}M_f'(\beta)=\psi(\beta)-\psi(\gamma-\beta).
$$
Hence by \eqref{e:derivative-at-zero},
\begin{equation}\label{e:derivative-at-zero-special}
i\overline{\Psi}'(0)=\Gamma(\alpha)\frac{\sin(\pi \alpha\widehat\rho)}{\pi}
\left(\pi\cot(\pi\alpha\widehat\rho)+
\psi(\beta)-\psi(\gamma-\beta)\right).
\end{equation}

\noindent
{\bf (a)} Suppose that $\gamma=1$ in \eqref{e:phi-special}. Then $\beta\in (0,1)$, and by the reflection formula for the digamma function, $\psi(1-\beta)-\psi(\beta)=\pi\cot(\pi \beta)$,  see \cite[8.365.8]{GR07}.  Elementary calculation gives that
\begin{equation}
i\overline{\Psi}'(0)=-\Gamma(\alpha)\frac{\sin(\pi(\alpha\wh{\rho}-\beta))}{\sin(\pi \beta)}, \nonumber
\end{equation}
and the sign of $i\overline{\Psi}'(0)$ depends on the sign of $\sin(\pi(\alpha\wh{\rho}-\beta))$.

\noindent
{\it Case 1}: $\alpha=1$. Then $\wh{\rho}=1/2$, and $\sin(\pi(\alpha\wh{\rho}-\beta))=\sin(\pi/2-\pi\beta)=\cos(\pi \beta)$. Therefore, $i\overline{\Psi}'(0)>0$ if $\beta\in (1/2,1)$, $i\overline{\Psi}'(0)=0$ if $\beta=1/2$, and $i\overline{\Psi}'(0)<0$ if $\beta\in (0,1/2)$. 

\noindent
{\it Case 2}: $\alpha\in (0,1)\cup (1,2)$. Since $\alpha\wh{\rho}\in (0,1)$, we have that $\alpha\wh{\rho}-\beta\in (-1,1)$. Thus, $\sin(\pi(\alpha\wh{\rho}-\beta))=0$ if and only if $\rho=1-\frac{\beta}{\alpha}$. Since we must have that $\rho>0$ and $\rho< 1/\alpha$, we get two critical values:  $\alpha_*=\beta$ and $\alpha^*:=1+\alpha_*=1+\beta$, see Corollary \ref{c:s2}. If $\alpha\in (0, \alpha_*)$, then $i\overline{\Psi}'(0)>0$; If $\alpha\in [\alpha_*, \alpha^*]$, then $i\overline{\Psi}'(0)>0$ if $\rho>1-\frac{\beta}{\alpha}$, $i\overline{\Psi}'(0)=0$ if $\rho=1-\frac{\beta}{\alpha}$, and $i\overline{\Psi}'(0)<0$ if $\rho<1-\frac{\beta}{\alpha}$; If $\alpha\in ( \alpha^*,2)$, then $i\overline{\Psi}'(0)<0$. 

\begin{figure}[h!]
	\begin{tabular}{cc}
        \includegraphics[width=8.0cm,height=5.0cm]{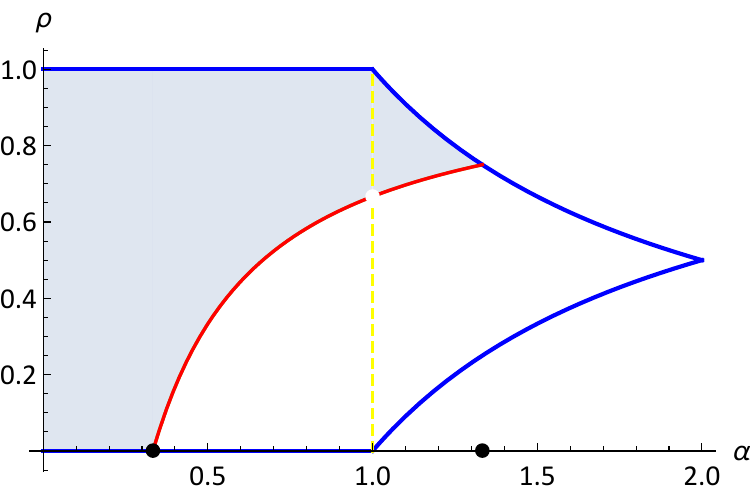}&\includegraphics[width=8.0cm,height=5.0cm]{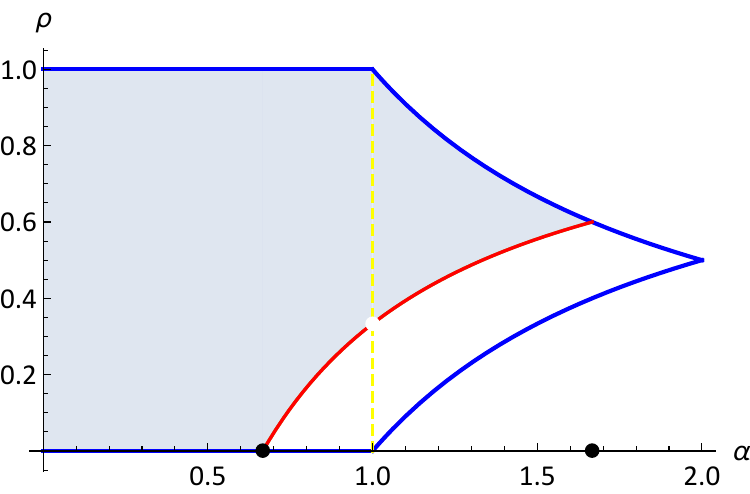}
	\end{tabular}
    \caption{Example \ref{ex:beta-gamma} (a): Left: $\beta=1/3$, $\gamma=1$, $\alpha_*=1/3$, $\alpha^*=4/3$; Right: $\beta=2/3$, $\gamma=1$, $\alpha_*=2/3$, $\alpha^*=5/3$;     $i\overline{\Psi}'(0)>0$ in the shaded region, $i\overline{\Psi}'(0)=0$ on the red line.  }
\end{figure}

We further assume that $\beta=1-\alpha\rho$. The pssMp $\overline{X}$ is then the part (until the first hitting of zero) of 
the trace process of the $\alpha$-stable process $\eta$ in $(0, \infty)$ (or the  path-censored $\alpha$-stable process). 
We have  that $\alpha \wh{\rho}-\beta =\alpha \wh{\rho}+\alpha \rho-1 =\alpha-1$. Thus, if $\alpha<1$, then $\alpha \wh{\rho}-\beta\in (-1,0)$, and therefore $i\overline{\Psi}'(0)>0$. If $\alpha=1$, then $\alpha \wh{\rho}-\beta=0$, and therefore $i\overline{\Psi}'(0)=0$. If $\alpha>1$, then $\alpha \wh{\rho}-\beta\in (0,1)$, and therefore $i\overline{\Psi}'(0)<0$. This shows that $\overline{X}$ has infinite 
absorption time in case $\alpha\in(0,1]$ and hits zero in finite time when $\alpha\in (1,2)$. Of course, since $\overline{X}$ is the trace process, this fact is well known. 

\noindent
{\bf (b)}
Let us now consider the case $\beta=1$ and $\gamma=\alpha+1$.
In this case
\begin{eqnarray}\label{e:Psi'-for-Z}
i\overline{\Psi}'(0)
&=&\frac{\Gamma(\alpha)\sin(\pi \alpha\widehat\rho)}{\pi}\left(\pi\cot(\pi \alpha\widehat\rho))+(\psi(1)-\psi(\alpha))\right).
\end{eqnarray}
The equation $\pi\cot(\pi(\alpha-1))=\psi(\alpha)-\psi(1)$ (obtained by formally taking $\rho=1/\alpha$), has a unique solution $\alpha^{\ast}$  in $(1,2)$ which can be numerically computed. It turns out that $\alpha^{\ast}\approx 1.44386$ with corresponding $\rho^{\ast}=1/\alpha^{\ast}\approx 0.692588$. Further, solving $\pi\cot(\pi \alpha (1-\rho))+(\psi(1)-\psi(\alpha))=0$ for $\rho$, we get a unique solution
$$
\rho(\alpha)=1-\frac{1}{\alpha \pi}\mathrm{arccot} \left(\frac{\psi(1)-\psi(\alpha)}{\pi}\right).
$$
It is easy to see that  $\lim_{\alpha \downarrow 0}\rho(\alpha)=-\infty$, $\lim_{\alpha\uparrow 2}\rho(\alpha)=1$, hence by continuity there exists $\alpha_*$ such that $\rho(\alpha_*)=0$, and consequently, $\rho(\alpha)>0$ for $\alpha\in (0,\alpha_*)$. Numerically we obtain $\alpha_*\approx 0.596051$. 

Therefore, we conclude that  (i) If $\alpha\in (0,\alpha_*]$, then $i\overline{\Psi}(0)>0$ for all $\rho\in (0,1)$;  (ii) if  $\alpha\in (\alpha_*, \alpha^*)$, then $i\overline{\Psi}(0)>0$ for $\rho>\rho(\alpha)$, $i\overline{\Psi}(0)=0$ for $\rho=\rho(\alpha)$, $i\overline{\Psi}(0)<0$ for $\rho<\rho(\alpha)$; (iii) if $\alpha\in [\alpha^*, 2)$, then $i\overline{\Psi}(0)<0$ for all admissible $\rho$. 
\begin{figure}[h!]
   \begin{center}
     \includegraphics[width=8.0cm,height=5.0cm]{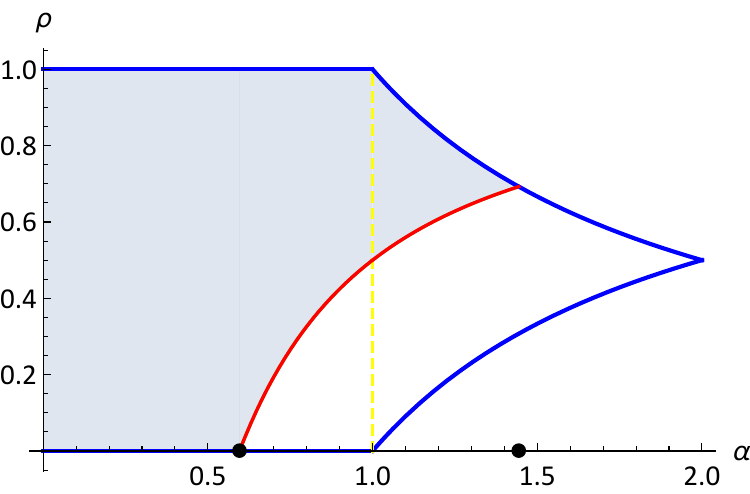}
    \end{center}
    \caption{Example 
    \ref{ex:beta-gamma} (b): $\alpha_*\approx 0.596501$, $\alpha^*\approx 1.44386$; $i\overline{\Psi}'(0)>0$ in the shaded region, $i\overline{\Psi}'(0)=0$ on the red line.}\label{cap:2}
\end{figure}
}
\end{example}

\begin{example}\label{ex:reflected-process}
{\rm
In case $\phi=\delta_a$, 
\begin{equation}\label{e:derivative-at-zero-case2}
i\overline{\Psi}'(0)=\Gamma(\alpha)\frac{\sin(\pi \alpha\widehat\rho)}{\pi}
\left(\pi\cot(\pi\alpha\widehat\rho)-\log a\right).
\end{equation}
If $a=1$, the corresponding pssMp $\overline{X}$ of index $\alpha$  is the absolute value of the strictly $\alpha$-stable L\'evy process. The sign of $i\overline{\Psi}'(0)$ depends on whether $\alpha\wh{\rho}$ is less than, equal, or larger than 1/2. Similar analysis can be done for any $a>0$. 
\begin{figure}[h!]
   \begin{center}
        \includegraphics[width=8.0cm,height=5.0cm]{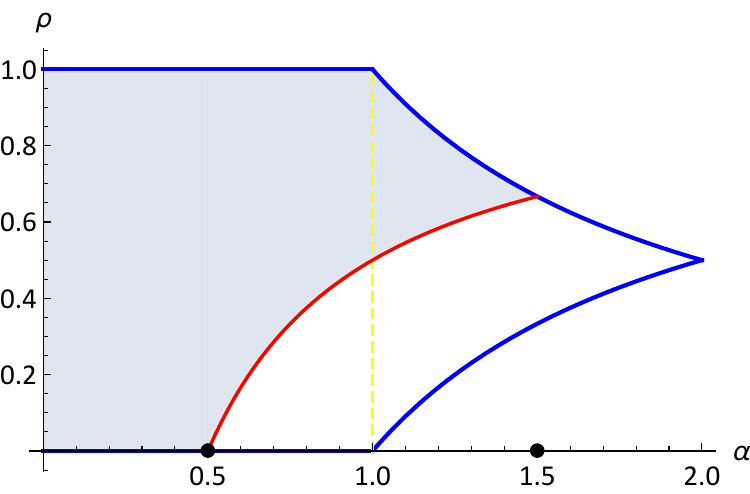}
    \end{center}
    \caption{Example \ref{ex:reflected-process}: $a=1$; $i\overline{\Psi}'(0)>0$ in the shaded region, $i\overline{\Psi}'(0)=0$ on the red line.}
\end{figure}
}
\end{example}

\begin{example}\label{ex:phi-exp-decay}{\rm
We now look at the example in which  $\phi$ is given by \eqref{e:phi-special-exp}.
Then $\phi$ is a probability density on $(0,{\infty})$. Its Mellin transform is given by
\begin{eqnarray*}
M_{\phi}(s)&:=&\int_0^{\infty}\phi(t) t^{s-1}dt=\frac{\gamma a^{\frac{\beta}{\gamma}}}{\Gamma(\frac{\beta}{\gamma})}\int_0^{\infty}t^{s+\beta-2} e^{-at^{\gamma}}dt\\
&=&\frac{\gamma a^{\frac{\beta}{\gamma}}}{\Gamma(\frac{\beta}{\gamma})}\, \gamma^{-1} a^{-\frac{s+\beta-1}{\gamma}} \Gamma\left(\frac{s+\beta-1}{\gamma}\right)
=\frac{ \Gamma\left(\frac{s+\beta-1}{\gamma}\right)}{\Gamma\left(\frac{\beta}{\gamma}\right)}\, a^{\frac{1-s}{\gamma}}.
\end{eqnarray*}
Thus we have
$$
\int_0^{\infty}\phi(t)(\log t)t^{s-1}dt=M_{\phi}'(1)=\frac{ \Gamma\left(\frac{s+\beta-1}{\gamma}\right)}{\gamma\Gamma\left(\frac{\beta}{\gamma}\right)}\left(\psi\left(\frac{s+\beta-1}{\gamma}\right)- \log a\right)a^{\frac{1-s}{\gamma}}, 
$$
and finally,
$$
\int_0^{\infty}\phi(t)(\log t)dt=M_{\phi}'(1)=\frac{1}{\gamma}\left(\psi\left(\frac{\beta}{\gamma}\right)- \log a\right).
$$
Therefore,
$$
i\overline{\Psi}'(0)=\Gamma(\alpha)\frac{\sin(\pi \alpha\widehat\rho)}{\pi}
\left(\pi\cot(\pi\alpha\widehat\rho)-\frac{1}{\gamma}\left(\psi\left(\frac{\beta}{\gamma}\right)- \log a\right)\right).
$$
\begin{figure}[h!]
   \begin{tabular}{cc}
      \includegraphics[width=8.0cm,height=5.0cm]{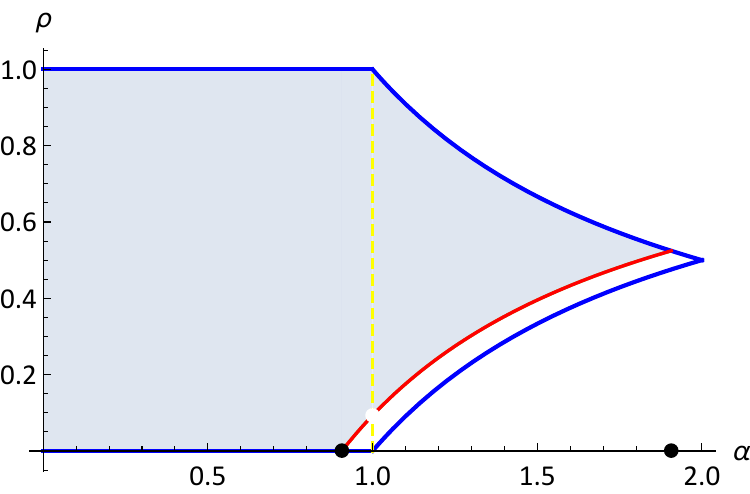}&\includegraphics[width=8.0cm,height=5.0cm]{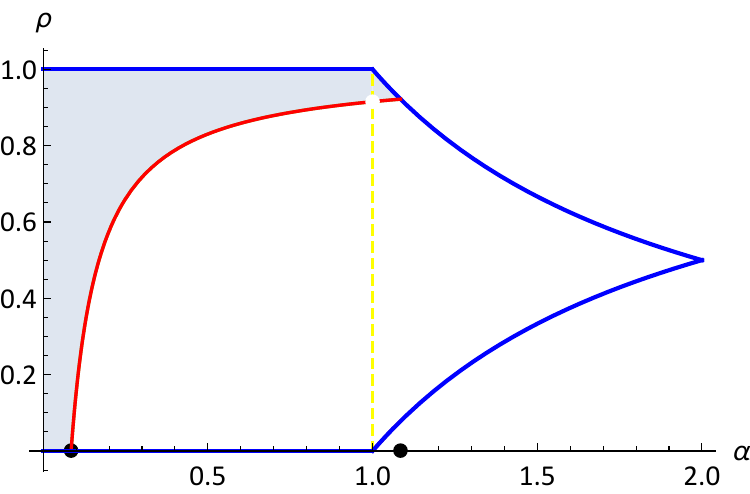}
      \end{tabular}
    \caption{Example \ref{ex:phi-exp-decay}: Left: $a=1$, $\beta=0.1$, $\gamma=1$, $\alpha_*\approx 0.906821$, $\alpha^*\approx 1.906821$; Right: $a=1$, $\beta=100000$, $\gamma=1$, $\alpha_*\approx 0.0847945$, $\alpha^*\approx 1.0847945$; $i\overline{\Psi}'(0)>0$ in the shaded region, $i\overline{\Psi}'(0)=0$ on the red line. }
\end{figure}
}
\end{example}

We end this subsection with the analysis of the behavior at 
lifetime 
(absorption time)
of the censored process. To the best of our knowledge, this has not been completely done before, but see \cite[p.976]{CC06}.

\begin{example}\label{ex:censored}
{\rm
In this example we consider the L\'evy process with characteristic exponent
$$
\Psi(\theta)=\frac{\Gamma(\alpha-i\theta)}{\Gamma(\alpha\wh{\rho}-i\theta)}\, \frac{\Gamma(i\theta+1)}{\Gamma(i\theta+1-\alpha\wh{\rho})}-\frac{c_-}{\alpha}.
$$
The corresponding pssMp $X$ is 
the (not necessarily symmetric) censored $\alpha$-stable process.
It is straightforward to calculate that
\begin{eqnarray}\label{e:Psi'-for-censored}
i\Psi'(0)
=\frac{\Gamma(\alpha)\sin(\pi \alpha\widehat\rho)}{\pi}\left(\pi\cot(\pi \alpha\widehat\rho)-(\psi(1)-\psi(\alpha))\right).
\end{eqnarray}
Notice the similarity with the expression in Example 
 \ref{ex:beta-gamma} (b) and
the difference being the sign in front of $(\psi(1)-\psi(\alpha))$. 
The equation $\pi\cot(\pi(\alpha-1))=\psi(1)-\psi(\alpha)$ (obtained by formally taking $\rho=1/\alpha$), has a unique solution $\alpha^{\ast}$  in $(1,2)$ which can be numerically computed. It turns out that $\alpha^{\ast}\approx 1.56735$ with corresponding $\rho^{\ast}=1/\alpha^{\ast}\approx 0.63802$. For every $\alpha\in (0,\alpha^{\ast})$ equation
$\pi \cot(\pi\alpha(1-\rho))=\psi(1)-\psi(\alpha)$ has a unique solution given by
$$
\rho(\alpha)=1-\frac{1}{\alpha \pi}\mathrm{arccot} \left(\frac{\psi(1)-\psi(\alpha)}{\pi}\right)
$$
which is strictly increasing in $\alpha$.  
Moreover, it can be shown that $\rho(\alpha)>0$ for every $\alpha\in (0,2)$, and $\lim_{\alpha\downarrow 0}\rho(\alpha)=0$ (so formally we can take $\alpha_*=0$). 
When $\rho>\rho(\alpha)$ we have $i\Psi'(0)>0$, for $\rho=\rho(\alpha)$ it holds that $i\Psi'(0)=0$, while for $\rho<\rho(\alpha)$, $i\Psi'(0)<0$. When $\alpha\in [\alpha^{\ast}, 2)$, for every admissible $\rho$ we have that $i\Psi'(0)<0$. Note also that for every $\rho\in [\rho^{\ast},1)$ it holds that $i\Psi'(0)>0$. 
\begin{figure}[h!]
   \begin{center}
     \includegraphics[width=8.0cm,height=5.0cm]{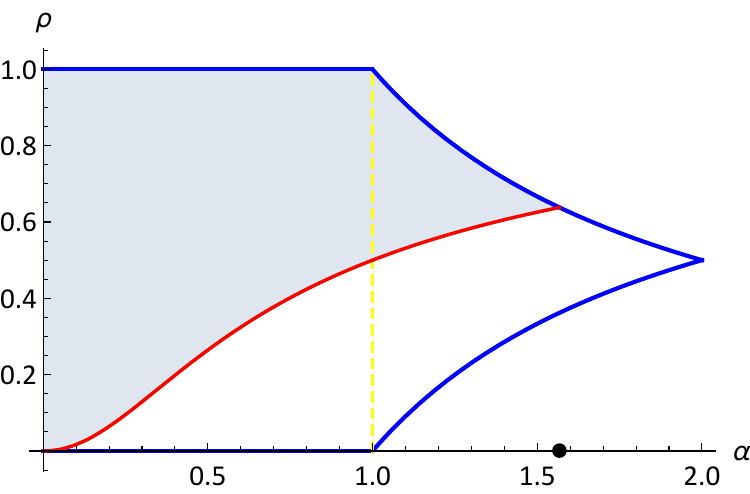}
    \end{center}
    \caption{Example \ref{ex:censored}: $\alpha^*\approx 1.56735$; $i\overline{\Psi}'(0)>0$ in the shaded region, $i\overline{\Psi}'(0)=0$ on the red line.}
\end{figure}
}
\end{example}


\section{Symmetric  resurrection kernels}\label{s:sym-int-kernel}
The goal of this section is to find a sufficient and necessary condition for the  resurrection kernel to be symmetric, namely that it holds $q(x,y)=q(y,x)$. 
When the  resurrection kernel is symmetric and the strictly $\alpha$-stable process $\eta$ is also symmetric (that is $c_+=c_-$, or, equivalently, $\rho=1/2$), 
the jump kernel $J(x,y)=j(x,y)+q(x,y)$ of the resurrected process $\overline{X}$ is also symmetric. In particular, the process $\overline{X}$ is symmetric with respect to the Lebesgue measure in $(0,{\infty})$.

We first give a necessary and sufficient condition for symmetry in terms of the L\'evy measure $\pi$ of the compound Poisson process $\chi$.

\begin{prop}\label{p:q-symmetric-levy}
Let $\pi$ be the L\'evy measure of the compound Poisson process $\chi$. Then $q(x,y)=q(y,x)$ for all $x,y>0$, if and only if,
\begin{equation}\label{e:q-symmetric-levy}
\pi(-y)=e^{(\alpha-1)y}\pi(y).
\end{equation}
\end{prop}
\pf Recall that $\pi(y)=e^{-\alpha y}q(e^{-y},1)$. Suppose that $q$ is symmetric. Then by  symmetry and 
scaling
\begin{eqnarray*}
\pi(-y)&=&e^{\alpha y}q(e^y,1)=e^{\alpha y}q(1,e^y)=e^{\alpha y} (e^y)^{-1-\alpha}q(e^{-y},1)\\
&=&e^{-y}q(e^{-y},1)=e^{-y} e^{\alpha y}\pi(y)=e^{(\alpha-1)y}\pi(y).
\end{eqnarray*}
The converse is similar. \qed

Recall from Lemma \ref{l:q-always-density} that
$$
q(x,y)=c_- 
\int_{(0,{\infty})} 
\left(x+\frac{y}{v}\right)^{-1-\alpha}\frac{\phi(dv)}{v}.
$$
 
The proof of the next technical lemma is given in the appendix.
\begin{lemma}\label{l:q-symmetric}
Suppose that $m$ is a signed Borel measure on $(0, \infty)$ such that
$$
\int_{(0, \infty)}(1+xu)^{-1-\alpha}|m|(du)<\infty, \quad  \text{for all }   x>0
$$
and
\begin{equation}\label{e:unique1}
\int_{(0, \infty)}(1+xu)^{-1-\alpha}m(du)=0, \quad  \text{for all }  x>0.
\end{equation}
Then $m$ is the zero measure on $(0, \infty)$.
\end{lemma}

\begin{thm}\label{t:q-symmetric}
It holds that $q(x,y)=q(y,x)$ for all $x,y>0$, if and only if
\begin{equation}\label{e:q-symmetric}
\phi_*(dt)=t^{\alpha-1}\phi(dt), \quad \mbox{ on }  (0, \infty),  
\end{equation}
where $\phi_*$ is the pushforward of 
the restriction of  the measure $\phi$  to $(0, \infty)$ 
under the map $x\to 1/x$. In case when 
the restriction of the measure $\phi(dt)$  to $(0, \infty)$
has a density $\phi(t)$ with respect
to the Lebesgue measure, the measure equality above reduces to $\phi(t^{-1})=t^{1+\alpha}\phi(t)$ for almost every $t>0$.
\end{thm}
\pf We have that
\begin{eqnarray*}
q(y,x)&=&c_- 
\int_{(0, \infty)}
\left(y+\frac{x}{v}\right)^{-1-\alpha}\frac{\phi(dv)}{v}=
c_- \int_{(0, \infty)}
(yv+x)^{-1-\alpha} v^{\alpha}\phi(dv) \\
&=&c_-\int_{(0, \infty)}
\left(\frac{y}{u}+x\right)^{-1-\alpha}u^{-\alpha}\phi_*(du).
\end{eqnarray*}
If \eqref{e:q-symmetric} holds, then the last integral in the display above is equal to
$$
c_- \int_{(0, \infty)}
\left(\frac{y}{u}+x\right)^{-1-\alpha}u^{-\alpha}u^{\alpha-1}\phi(du) =c_- \int_{(0, \infty)}
\left(x+\frac{y}{u}\right)^{-1-\alpha}\frac{\phi(du)}{u} =q(x,y).
$$

Conversely, assume that $q$ is symmetric. Then we must have that 
$$
\int_{(0, \infty)}
\left(\frac{y}{u}+x\right)^{-1-\alpha} u^{-\alpha}\phi_*(du)= 
\int_{(0, \infty)}
\left(\frac{y}{u}+x\right)^{-1-\alpha}\frac{\phi(du)}{u},
$$
for all $x,y>0$. By taking $y=1$ and rewriting, we get that
$$
\int_{(0, \infty)}
(1+xu)^{-1-\alpha}u\phi_*(du)= 
\int_{(0, \infty)}
(1+xu)^{-1-\alpha} u^{\alpha}\phi(du),\quad \text{for all }x>0.
$$
The claim now follows from Lemma \ref{l:q-symmetric}. \qed

\begin{corollary}\label{c:phi-sym-general}
Let $\phi:(0, \infty)\to [0,\infty)$ be such that $\int_0^{\infty}\phi(t)dt=1$. Then $\phi$ satisfies \eqref{e:q-symmetric} if and only if
\begin{equation}\label{e:phi-sym-general}
\phi(t)=\frac{f(t+t^{-1})}{(1+t)^{1+\alpha}},
\end{equation}
for $f:[2,\infty)\to[0,\infty)$. 
\end{corollary}
\pf Assume that $\phi$ is given by \eqref{e:phi-sym-general}. Let $g(t):=f(t+t^{-1})$. Then $g(t^{-1})=g(t)$ implying that $\phi$ satisfies  \eqref{e:q-symmetric}. 

Conversely, if $\phi$ satisfies  \eqref{e:q-symmetric}, define 
$g(t):=\phi(t)(1+t)^{1+\alpha}$.
Then $g(t^{-1})=g(t)$. 
For $s\ge 2$, solving $t+t^{-1}=s$, we
 get two solutions: $t=(s+\sqrt{s^2-4})/2 \ge1$ and $t^{-1}=(s-\sqrt{s^2-4})/2\le 1$. Define $f(s):=g(t)=g(t^{-1})$. Since $s=t+t^{-1}$, we see that $g(t)=f(t+t^{-1})$. \qed

\begin{prop}\label{p:region-symmetry} 
Suppose that $\mathfrak{p}=1$ and  $q$ is symmetric. 
Then  ${\bf E}[\overline{\xi}_1]< 0$ if $\alpha > 1$  
and $\widehat\rho \in  [1/(2\alpha), 1/\alpha)$, 
${\bf E}[\overline{\xi}_1]=0$ if $\alpha =1$ and $\widehat\rho = 1/2$, and 
${\bf E}[\overline{\xi}_1]>0$ if $\widehat\rho \in  (0,  1/(2\alpha)]$.

In particular, if $\alpha \ge 3/2$, then the 
absorption time of $\overline{X}$ is finite $\P_x$-a.s. and $X$ is continuously 
absorbed at 0. Also, if $\alpha\le 1/2$, then the
absorption time of $\overline{X}$ is infinite.
\end{prop}
\pf
By Theorem \ref{t:q-symmetric},
\begin{align*}
\int_{(0, \infty)}
(\log u)\phi(du)&=
\int_{(0, 1)}
(\log u)\phi(du)+
\int_{(1, \infty)}
(\log u)\phi(du)\\
&=\int_{(1, \infty)}(\log 1/u)\phi_*(du)+\int_{(1, \infty)}(\log u)\phi(du) \\
&=\int_{(1, \infty)}
(\log u)(-u^{\alpha-1}+1)\phi(du),
\end{align*}
which is negative  for  $\alpha > 1$, zero for   $\alpha =1$, and positive for $\alpha < 1$.  It follows from \eqref{e:derivative-at-zero} that the sign of $i\overline{\Psi}'(0) $ depends on the sign of 
$$
\pi\cot(\pi\alpha  \widehat \rho)-
\int_{(1, \infty)}
(\log u)(u^{\alpha-1}-1)\phi(du).
$$
This expression is negative for $\alpha > 1$ and $\widehat\rho \in  [1/(2\alpha), 1/\alpha)$,  zero for $\alpha =1$ and $\widehat\rho = 1/2$,
and positive for $\alpha <1$ and $\widehat\rho \in  (0,  1/(2\alpha)]$.

The absorption claim for $\alpha\ge 3/2$ (which is equivalent to $1-1/\alpha\ge 1/(2\alpha)$) follows from the assumption $\wh{\rho}>1-1/\alpha$.
\qed

\begin{figure}[h!]
   \begin{center}
     \includegraphics[width=8.0cm,height=5.0cm]{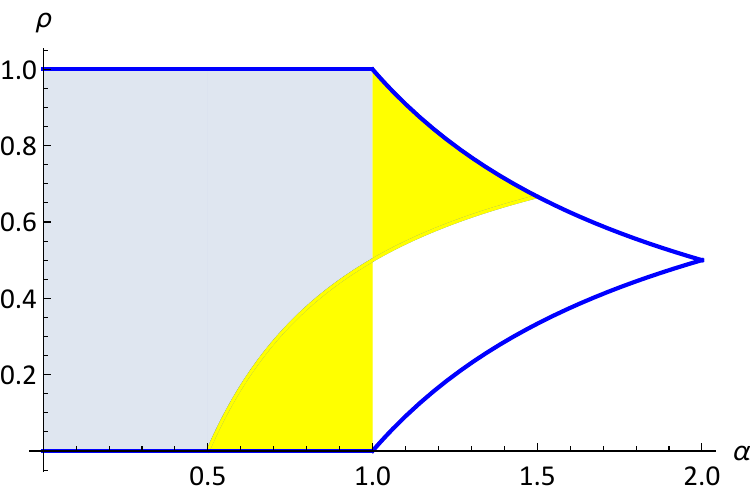}
    \end{center}
    \caption{Illustration of Proposition \ref{p:region-symmetry}: In the shaded region $i\overline{\Psi}'(0)>0$, in the white region  
     $i\overline{\Psi}'(0) <  0$, in the yellow region the sign of $i\overline{\Psi}'(0)$ is undetermined.}
\end{figure}

\begin{example}\label{ex:symmetric-phi}
{\rm
(a) Let $\phi$ be as in \eqref{e:phi-special}:
$$
\phi(t)=\frac{\Gamma(\gamma)}{\Gamma(\beta)\Gamma(\gamma-\beta)}t^{\beta-1}(1+t)^{-\gamma}.
$$
Imposing condition \eqref{e:q-symmetric} on $\phi$ implies that $\gamma=\alpha+2\beta-1$. Since $\gamma>\beta$, we get that $\beta >1-\alpha$. Thus
$$
\phi(t)=\frac{\Gamma(\alpha+2\beta-1)}{\Gamma(\beta)\Gamma(\alpha+\beta-1)}   t^{\beta-1}(1+t)^{1-\alpha-2\beta}.
$$

Note that if $\beta=1$ then symmetry requires that $\gamma=\alpha+1$, so we get Example 
 \ref{ex:beta-gamma} (b). In this case, the regions for the sign of  $i\overline{\Psi}'(0)$ are
completely determined, see Figure \ref{cap:2}. 

\noindent
(b) If $\phi$ is as in \eqref{e:phi-special-exp}, then it cannot lead to a symmetric  resurrection kernel.

\noindent 
(c) Let $\phi=\delta_a$, $a>0$. Then the measure $\phi$ satisfies \eqref{e:q-symmetric} if and only if $a=1$. 
}
\end{example}


\section{Sharp two-sided  estimates of the  resurrection kernel }\label{s:estimates-q}
In this section we establish sharp two-sided estimates of $q(x,y)$ under minimal assumptions. First notice that it follows from \eqref{e:q-always-density} that $q$ enjoys the following scaling
$$
q(\lambda x, \lambda y)=\lambda^{-1-\alpha}q(x,y), \quad x,y >0, \ \lambda>0. 
$$
This implies that for all $0<x<y$ we have
\begin{align}
q(x,y)=&(y-x)^{-1-\alpha}q\left(\frac{x}{y-x}, \frac{x}{y-x}+1\right) \label{e:scalq1}\\
q(y,x)=&(y-x)^{-1-\alpha}q\left(\frac{x}{y-x}+1, \frac{x}{y-x}\right). \label{e:scalq2}
\end{align}
Thus it suffices to get the estimates of $q(x,x+1)$ and $q(x+1,x)$, $x>0$.

We first look at the simple case when the measure $\phi$ has compact support in $(0,{\infty})$. Then it is easy to see that 
$$
q(x,x+1)=c_- \int_{(0, \infty)}
\left(x+1+\frac{x}{t}\right)^{-1-\alpha}  \frac{\phi(dt)}{t}\asymp 1
$$
and 
$$
q(x+1,x)=c_-\int_{(0, \infty)}
\left(x+\frac{x+1}{t}\right)^{-1-\alpha}  \frac{\phi(dt)}{t}\asymp 1.
$$
Thus in this case by \eqref{e:scalq1}-\eqref{e:scalq2}, we have
$$
q(x, y)\asymp |x-y|^{-1-\alpha}, \quad x, y\in (0, \infty).
$$

In the remainder of this section, we assume that $\phi$ is absolutely continuous
and has a strictly positive density.

\medskip
Assume that $|y-x|=1$ so that either $y=x+1$ or $y=x-1$. Then by Lemma \ref{l:q-always-density},
$$
q(x+1,x)
= c_- \int_0^{\infty}\left(x+1+\frac{x}{t}\right)^{-1-\alpha}  \phi(t)\frac{dt}{t}
$$
and 
$$
q(x,x+1)
= c_-\int_0^{\infty}\left(x+\frac{x+1}{t}\right)^{-1-\alpha}  \phi(t)\frac{dt}{t}.
$$

\begin{lemma}\label{l:estimates-of-q}
(1) If $x\ge 1/4$, then
$
q(x+1,x)\asymp q(x,x+1)\asymp x^{-1-\alpha}.
$

\noindent
(2) Suppose $\phi$ satisfies the 
lower weak scaling condition 
$L_1(\beta_1)$ 
at zero with $\beta_1>-1-\alpha $.
Then for $x\le 1/4$, 
\begin{equation}\label{e:estimate-q-1}
q(x+1,x)\asymp \int_{x}^{1}  \phi(t)\frac{dt}{t}.\end{equation}
Further, 
if  
$\phi$ also satisfies the upper weak scaling condition 
$U_1(\beta_2)$ 
at zero with $\beta_2<0$, then 
\begin{equation}\label{e:estimate-q-10}
q(x+1,x)\asymp  \phi(x).\end{equation}
(3) Suppose $\phi$ satisfies the upper weak scaling condition  
$U^1(\gamma_2)$ 
at infinity  with  $\gamma_2<0$. Then
\begin{equation}\label{e:estimate-q-2}
q(x,x+1)
\asymp \int_0^{1/x} t^{\alpha} \phi(t){dt} \asymp \int_1^{1/x} t^{\alpha} \phi(t){dt}.\end{equation}
Further, if  $\phi$ also satisfies the lower weak scaling condition 
$L^1(\gamma_1)$ 
at infinity  with $\gamma_1<-1-\alpha$,
\begin{equation}\label{e:estimate-q-20}
q(x,x+1)
\asymp x^{-1-\alpha} \phi(1/x).\end{equation}
\end{lemma}
\pf Without loss of generality, we neglect the constant $c_-$ in the proof.
First note that 
\begin{align}\label{e:est1}
0<\int_0^{\infty}\frac{t^\alpha\phi(t)}{(1+t)^{\alpha+1}}dt  \le 
\int_0^{\infty}\phi(t)dt =1.
\end{align}

\noindent (1)
{\bf Case 1:} $1/4 \le x \le 4$.
Then 
$$
x+1+\frac{x}{t}\ge 1+\frac{1}{4t}\ge \frac14\left(1+\frac{1}{t}\right) \quad \text{ and} \quad x+1+\frac{x}{t} \le 4+1 +\frac{1}{t}=5\left(1+\frac{1}{t}\right).
$$
  Also,   
$$
x+\frac{x+1}{t}\ge \frac14 +\frac{1}{t} \ge \frac14 \left(1+\frac{1}{t}\right) \quad \text{ and} \quad  x+\frac{x+1}{t} \le 4 +\frac{4+1}{t}\le 5\left(1+\frac{1}{t}\right).
$$ 
Thus, by 
\eqref{e:est1}
we get
$$
q(x+1,x)\asymp q(x,x+1) \asymp \int_0^{\infty}\left(1+\frac{1}{t}\right)^{-1-\alpha}\phi(t)\frac{dt}{t} =
\int_0^{\infty}\frac{t^\alpha\phi(t)}{(1+t)^{\alpha+1}}dt \asymp 1.
$$

\medskip
\noindent {\bf Case 2:} $ x \ge 4$.
Then we have that  
$$
x+\frac{x}{t}\le x+1+\frac{x}{t}\le 2x+2\frac{x}{t},
$$
and hence 
$$
x+1+\frac{x}{t} \asymp x\left(1+\frac{1}{t}\right). 
$$
Also,  
$$
x+\frac{x+1}{t}\le (x+1)\left(1+\frac{1}{t}\right) \quad \text{ and } \quad x+\frac{x+1}{t}\ge 
\frac{x+1}{2}\left(1+\frac{1}{t}\right)
$$ 
imply 
$$
x+\frac{x+1}{t}\asymp (x+1)\left(1+\frac{1}{t}\right)\asymp x\left(1+\frac{1}{t}\right).
$$
Thus, by \eqref{e:est1}
$$
q(x+1,x)\asymp q(x,x+1) \asymp x^{-1-\alpha} \int_0^{\infty}\left(1+\frac{1}{t}\right)^{-1-\alpha}\phi(t)\frac{dt}{t}\asymp  x^{-1-\alpha}.
$$

\medskip
\noindent
(2) Assume $x \le 1/4$ and let 
$$
q(x+1,x)=\int_0^{x}+ \int_{x}^{1}+\int_{1}^{\infty}=:I+II+III.
$$

For $0<t<x$, we have that 
$$
x+1=\frac{t(x+1)}{x}\frac{x}{t} \le (x+1)\frac{x}{t}\le 2\frac{x}{t},
$$
 hence $x+1+\frac{x}{t}\asymp \frac{x}{t}$. 
Since $\phi$ satisfies the lower weak scaling condition at zero
$L_1(\beta_1)$ 
with $\beta_1>-1-\alpha $,
we have
$$\int_0^{x} t^{\alpha} \phi(t){dt} =\phi(x)\int_0^{x} t^{\alpha} (\phi(t)/\phi(x)){dt}   \le c^{-1} 
\phi(x)\int_0^{x} t^{\alpha} (t/x)^{\beta_1}{dt} 
=c_1 x^{\alpha+1} \phi(x).$$
Thus,
\begin{align*}
I \asymp \int_0^{x} \left(\frac{x}{t}\right)^{-1-\alpha} \phi(t)\frac{dt}{t}
= x^{-1-\alpha}\int_0^{x} t^{\alpha} \phi(t){dt} 
\le c_1  \phi(x).
\end{align*}

For $x\le t <\infty$, we have 
$$
1\le x+1+\frac{x}{t}\le \frac14+1+1,
$$
hence $x+1+x/t\asymp 1$. 
Therefore,
\begin{eqnarray*}
II &\asymp & \int_{x}^{1}  \phi(t)\frac{dt}{t}. \end{eqnarray*}
Finally,
$$
III  \asymp \int_1^{\infty}\phi(t)\frac{dt}{t} \le  \int_1^{\infty}\phi(t){dt}  \le 1.
$$
Since
$$
\int_{x}^{1}  \phi(t)\frac{dt}{t} 
\ge \phi(x) \int_{x}^{2x} \frac{\phi(t)dt}{ \phi(x)t}+  2 \int_{1/2}^{1} \phi(t){dt}
\ge  c 
\phi(x) \int_{x}^{2x} \frac{t^{\beta_1-1}}{x^{\beta_1}}dt+c_2
 \asymp  \phi(x)+1,
 $$ 
we get \eqref{e:estimate-q-1}. Moreover,  if  $\phi$ also satisfies 
$U_1(\beta_2)$ 
with $\beta_2<0$, 
then 
$$
c_3 \phi(x) = c \phi(x)\int_{x}^{1}  (t/x)^{\beta_1}\frac{dt}{t} \le
\int_{x}^{1}  \phi(t)\frac{dt}{t} \le  C  \phi(x)\int_{x}^{1}  (t/x)^{\beta_2}\frac{dt}{t} 
=c_4  \phi(x),
$$ 
so we get \eqref{e:estimate-q-10}.

\medskip
(3) Assume $x \le 1/4$.
$$
q(x,x+1)= \int_0^{1/x}+\int_{1/x}^{\infty}=:I+II.
$$

For $0<t<2$, 
$$
x+\frac{x+1}{t}\le x+1 +\frac{x+1}{t} \le 2\frac{x+1}{t}+(x+1)t,
$$
hence 
$$
x+\frac{x+1}{t}\asymp \frac{x+1}{t}\asymp \frac{1}{t}.
$$
For $2<t <1/x$,  
$$
x+\frac{x+1}{t} \le \frac{1}{t}+\frac{2}{t}=\frac{3}{t},
$$
 hence $x+(x+1)/t\asymp 1/t$. Therefore,
$$
I \asymp \int_0^{1/x} \left(\frac{1}{t}\right)^{-1-\alpha} \phi(t)\frac{dt}{t} = \int_0^{1/x} t^{\alpha} \phi(t){dt}.
$$

When $t>1/x$, 
$$
\frac{x+1}{t}=\frac{x+1}{x}\, \frac{x}{t}\le \frac{2}{x}\, \frac{x}{t}\le (2t)\frac{x}{t}=2x.
$$
Thus, $x+(x+1)/t\le 3x$, hence $x+(x+1)/t \asymp x$. Moreover, 
using that fact that $\phi$ satisfies the upper weak scaling condition at infinity
$U^1(\gamma_2)$ 
with  $\gamma_2<0 $,  we have
$$
\int_{1/x}^{\infty} \phi(t)\frac{dt}{t}
=\phi(1/x)\int_{1/x}^{\infty} \frac{\phi(t)}{\phi(1/x)t}  dt
  \le C 
\phi(1/x)\int_{1/x}^{\infty} (tx)^{\gamma_2}\frac{dt}{t}
=c_5 \phi(1/x).$$ 
Therefore,
$$
II\asymp x^{-1-\alpha}\int_{1/x}^{\infty} \phi(t)\frac{dt}{t}  \le c_5 x^{-\alpha-1}\phi(1/x).
$$
Now, using 
\begin{align*}
\int_0^{1/x} t^{\alpha} \phi(t){dt} &\ge 
\phi(1/x)\int_{1/(2x)}^{1/x}\frac{\phi(t)}{\phi(1/x)} t^{\alpha} {dt}+\int_{0}^{1}  t^{\alpha} \phi(t){dt}\\
&\ge C^{-1}
\phi(1/x)\int_{1/(2x)}^{1/x}(xt)^{\gamma_2} t^{\alpha} {dt}+c_6\asymp  x^{-\alpha-1}\phi(1/x)+1,\end{align*}
we get \eqref{e:estimate-q-2}.
Finally,   if  $\phi$ also satisfies 
$L^1(\gamma_1)$ 
with $\gamma_1<-1-\alpha$, then
\begin{align*}
&c_7 x^{-\alpha-1}\phi(1/x) \le c 
\phi(1/x)\int_1^{1/x} t^{\alpha} (tx)^{\gamma_1} {dt}  \le 
\int_1^{1/x} t^{\alpha} \phi(t){dt} \\
& \le C 
\phi(1/x)\int_1^{1/x} t^{\alpha} (tx)^{\gamma_2} {dt} 
\le c_8 x^{-\alpha-1}\phi(1/x),
\end{align*}
and so
we get \eqref{e:estimate-q-20}.
\qed

\medskip
\noindent
\textbf{Proof of Theorem \ref{t:estimates-of-q}:}
The result follows immediately from Lemma \ref{l:estimates-of-q} and the scaling relations \eqref{e:scalq1}--\eqref{e:scalq2}. \qed

\medskip
As a consequence of Theorem \ref{t:estimates-of-q} we can get estimates of the L\'evy density $\pi(u)=q(1,e^u)e^u$. We state the next corollary in its simple form.
\begin{corollary}\label{c:estimates-of-pi}
Suppose that  $\phi$ satisfies both the lower and upper scaling conditions at zero and infinity as in Theorem \ref{t:estimates-of-q}. Then
$$
\pi(u)\asymp \left\{ \begin{array}{ll}
1,& |u| \le \log 5, \\
e^u\phi(e^u), & |u|>\log 5.
\end{array} \right.
$$
\end{corollary}

\medskip
We first apply Theorem \ref{t:estimates-of-q} to a generalization of the function $\phi$ given in \ref{e:phi-special}.
For $\beta \ge 0$, $\gamma \ge \beta$ and $\delta_+, \delta_- \in \R$, let
\begin{align}
\label{e:phasy}
\phi(t)\asymp \left(\log(e+t)\right)^{\delta_+}\left(\log(e+t^{-1})\right)^{\delta_-}{t^{\beta-1}(1+t)^{-\gamma}},
\end{align}
and $\int_0^{\infty}\phi(t)dt=1$, which implies that 
$\delta_-<-1$ if
 $\beta = 0$ and $\delta_+<-1$ if  $\gamma = \beta$.
Then 
$$
q(x,y)\asymp\int_0^{\infty}\left(x+\frac{y}{t}\right)^{-1-\alpha} \frac{\left(\log(e+t)\right)^{\delta_+}\left(\log(e+t^{-1})\right)^{\delta_-}}{t^{2-\beta}(1+t)^{\gamma}}\, dt.
$$
It is straightforward to check that the function $\phi$ satisfies the scaling conditions assumed in Theorem \ref{t:estimates-of-q}. Hence we have the following result.

\begin{corollary}\label{c:estimates-of-q}
Suppose $\phi$ is a probability density on $(0, \infty)$ satisfying   \eqref{e:phasy}.
If $x \le y\le 5x$, then
$$
q(x,y)\asymp q(y,x)\asymp x^{-1-\alpha}\asymp y^{-1-\alpha}.
$$
If $5x\le y$, then
\begin{equation}\label{e:estimate-q-3}
q(y,x)\asymp 
 y^{-1-\alpha}
\left\{ \begin{array}{ll}
(y/x)^{1-\beta}(\log (y/x))^{\delta_-}, & \beta<1,\\
\left\{\begin{array}{ll}
(\log (y/x))^{1+\delta_-}, & \delta_- > -1, \\
\log(\log( y/x)), & \delta_- =-1, \\
1, &\delta_-<-1,
\end{array}\right.
& \beta=1,\\
1, & \beta>1,
\end{array}
\right. 
\end{equation}
and
\begin{equation}
q(x,y)
\asymp 
y^{-1-\alpha}
\left\{ \begin{array}{ll}
(y/x)^{\alpha+\beta-\gamma}(\log (y/x))^{\delta_+}, &\gamma-\beta <\alpha ,\\
\left\{\begin{array}{ll}
(\log (y/x))^{1+\delta_+}, & \delta_+ > -1, \\
\log(\log(y/x)) , & \delta_+ =-1, \\
1, &\delta_+ <-1,
\end{array}\right.
&  \gamma-\beta =\alpha ,\\
1, & \gamma-\beta >\alpha .
\end{array}
\right. 
\end{equation}\label{e:estimate-q-4}
\end{corollary}
Note that in the case 
$\gamma-\alpha <\beta<1$, or in the case $\gamma-\alpha =\beta=1$ and $\delta_+ \wedge \delta_-\ge -1$,
both $q(x,y)$ and $q(y,x)$ explode when $x\to 0$. This leads to the following path interpretation: The intensity of jumps to and away from points near 0 is much higher than in case of the stable process. Thus, on average, large jump to and away from 
points near 0 are more probable.

In the symmetric case,
we have $\gamma=\alpha+2\beta-1$ and $\delta_+=\delta_-$, see Example \ref{ex:symmetric-phi} (a).
Hence $-\alpha-\beta+\gamma=-(1-\beta)$, and the estimates for $q(x,y)$ and $q(y,x)$ coincide.

\medskip
Now we  apply Theorem \ref{t:estimates-of-q} to the function $\phi$ given in \eqref{e:phi-special-exp}.

\begin{corollary}\label{c:estimates-of-qq}
Suppose $\phi$ is a probability density on $(0, \infty)$ satisfying 
\begin{align}
\label{e:phi_ex}
\phi(t) \asymp t^{\beta-1} e^{-at^{\gamma}}, \quad t>0
\end{align}
where $a,\gamma>0$ and $\beta>0$.
If $x \le y\le 5x$, then
$$
q(x,y)\asymp q(y,x)\asymp x^{-1-\alpha}\asymp y^{-1-\alpha}.
$$
If $5x\le y$, then
\begin{equation}\label{e:estimate-q-14a}
q(y,x)\asymp y^{-1-\alpha}
\left\{ \begin{array}{ll}
(y/x)^{-\beta+1}, & 0<\beta<1,\\
\log (y/x)& \beta=1,\\
1, & \beta>1,
\end{array}
\right. 
\end{equation}
and
\begin{equation}\label{e:estimate-q-14}
q(x,y)
\asymp y^{-1-\alpha}\int_1^{\frac{y}{x}} t^{\alpha+\beta-1} e^{-at^{\gamma}}{dt}
\asymp y^{-1-\alpha}.
\end{equation}
\end{corollary}
\pf
For the second comparison in \eqref{e:estimate-q-14}, 
see \eqref{e:estimate-q-21} and 
$$
0<
\int_1^{4} t^{\alpha+\beta-1} e^{-at^{\gamma}}{dt} \le \int_1^{\frac{y}{x}-1} t^{\alpha+\beta-1} e^{-at^{\gamma}}{dt} \le\int_1^{\infty} t^{\alpha+\beta-1} e^{-at^{\gamma}}{dt} <\infty, \quad
5x\le y.
$$
\qed

In case $\beta \le 1$, we see that $q(y,x)$ explodes as $x\to 0$, but $q(x,y)$ stays bounded. This means that the process will have tendency for big jumps 
to points close to the origin.

As a consequence 
of Corollary \ref{c:estimates-of-qq} 
we can derive that,  when a probability density $\phi$  on $(0, \infty)$ satisfies \eqref{e:phi_ex}, $\pi(u)\asymp 1$ for $|u|\le \log 5$, $\pi(u)\asymp e^{-u\alpha}$ for $u>\log 5$, and
$$
\pi(u)\asymp \left\{\begin{array}{ll}
e^{u\beta}, & 0<\beta <1, \\
|u|e^u, & \beta=1, \\
e^u, & \beta >1,
\end{array}\right.
$$
when $u<-\log 5$.


\section{Modified jump kernel}\label{s:modified-jump}
\subsection{General case}\label{ss:2-general}
Stable process 
conditioned to stay positive is a pssMp that can be regarded as a resurrected stable process. 
However, it does not fall
into the framework of resurrected stable processes of this paper. In this section we introduce a larger class of pssMps by modifying
the jump kernel of the pssMp $X=(X_t,\P_x)$ of index $\alpha$ defined 
in Subsection \ref{ss:sspcp}.
Thus, $X$ is a not necessarily symmetric censored process. Let $j(x,y)=\nu(y-x)$, where $\nu$ is defined in \eqref{e:stable-density}. We define a new jump kernel by
$$
J(x,y):=\sB(x,y)j(x,y),\quad x,y>0,
$$
where $\sB:(0,\infty)\times (0,\infty)\to (0,\infty)$ is a function satisfying the following properties: 

\noindent
\textbf{(B1)} Homogeneity: $\sB(\lambda x, \lambda y)=\sB(x,y)$ for all $x,y>0$ and all $\lambda >0$. 

\noindent
\textbf{(B2)} Integrability: (a) $y\mapsto e^{-\alpha y} \sB(1,e^y)$ is integrable at $\infty$ and $y\mapsto e^y \sB(1,e^y)$ is integrable at $-\infty$;  (b) 
$y\mapsto \sB(1,e^y)|y|^{1-\alpha}$ is integrable at 0; (c) For all $x \in (0, \infty)$, $y\mapsto \sB(x,y)$ is locally integrable in $(0, \infty)\setminus \{x\}$. 

\noindent
\textbf{(B3)} 
Regularity:
If $\alpha\in [1,2)$, there exist $\theta>\alpha-1$ and $C>0$ such that 
$$
|\sB(x,x)-\sB(x,y)|\le 
C\left(\frac{|x-y|}{x\wedge y}\right)^{\theta}
$$ 
for $|x-y| \le  (x\wedge y)/4$. 
If $\alpha<1$, there exists $C>0$ such that 
$
\sB(x,y)\le C$
for $|x-y| \le  (x\wedge y)/4$. 

Without loss of generality, from now on, we assume that $\sB(1,1)=1$.

For any $\sB(\cdot, \cdot)$ satisfying \textbf{(B1)}--\textbf{(B3)}, we will construct a pssMp with the jump kernel $J$ above via the Lamperti transform of a
certain L\'evy process.

\medskip
We first show that the jump kernel $J(x,y)=j(x,y)+q(x,y)$ of the resurrected process $\overline{X}$ is of the form introduced above. Indeed,  $J(x,y)$ can be rewritten as
$$
J(x,y)=j(x,y)+q(x,y)=j(x,y)\left(1+\frac{q(x,y)}{j(x,y)}\right)=\sB(x,y)j(x,y),
$$
where we define $\sB(x,y):=1+q(x,y)/j(x,y)$ for $y\neq x$, and $\sB(x,x)=1$. Clearly, $\sB(x,y)$ satisfies \textbf{(B1)}. Next,
by \eqref{e:levy-xi*} and \eqref{e:densities-Q-q},
$$
\sB(1, e^y)=1+\frac{q(1,e^y)}{j(1,e^y)}=1+\frac{\pi(y)}{e^y \nu(e^y-1)}=1+\frac{\pi(y)}{\mu(y)},
$$
so that $\sB(1,e^y)\mu(y)=\mu(y)+\pi(y)$ is a 
L\'evy density, 
i.e., $\int_{\R}(1\wedge y^2)\sB(1,e^y)\mu(y)dy<\infty$. 
Indeed, since $\mu(y)\asymp e^{-\alpha y}$ at $+\infty$, $\mu(y)\asymp e^{-y}$ at $-\infty$, and $\mu(y)\asymp |y|^{-1-\alpha}$ near zero, we know that \textbf{(B2)} holds.
Finally, it follows from Theorem \ref{t:estimates-of-q} that for $x<y<(5/4)x$ or $y<x<(5/4)y$ it holds that $q(x,y)\asymp q(y,x)\asymp x^{-1-\alpha}\asymp y^{-1-\alpha}$. Hence, if $|x-y|<(x\wedge y)/4$, 
\begin{eqnarray*}
|\sB(x,y)-\sB(x,x)|&=&\frac{q(x,y)}{j(x,y)}\le (c_+\vee c_-)|x-y|^{1+\alpha}q(x,y)\\
&\le & C|x-y|^{1+\alpha}(x^{-1-\alpha}\vee y^{-1-\alpha})=C\left(\frac{|y-x|}{x\wedge y}\right)^{1+\alpha}.
\end{eqnarray*}
Thus \textbf{(B3)} holds with $\theta=1+\alpha$.  

\medskip
As examples of 
this general setting we also mention the $\alpha$-stable process conditioned to stay positive and the $\alpha$-stable process conditioned to hit 0 continuously,
see \cite{Cha96, CD05}.
The jump kernel of the former is
$$
J(x, y)=\frac{y^{\alpha\widehat\rho}}{x^{\alpha\widehat\rho}}\, j(x, y), \quad x, y>0
$$
and the latter
$$
J(x, y)=\frac{y^{\alpha\widehat\rho-1}}{x^{\alpha\widehat\rho-1}}\, j(x, y), \quad x, y>0.
$$
It is straightforward to show that for 
every $\gamma\in (-1,\alpha)$, the function $\sB(x,y)=(y/x)^{\gamma}$ satisfies conditions \textbf{(B1)}-\textbf{(B3)}.  
In fact, 
\textbf{(B1)} and \textbf{(B2)}(b)--(c) clearly hold. 
For \textbf{(B2)}(a), $\sB(1, e^y)=e^{\gamma y}$, so $e^{-\alpha y}\sB(1,e^y)=e^{-(\alpha-\gamma)y}$ and is integrable at $\infty$ if and only if $\gamma <\alpha$. Also, $e^y \sB(1,e^y)=e^{(1+\gamma)y}$ and is integrable at $–\infty$ if and only if  
$\gamma >-1$.
For \textbf{(B3)}, without loss of generality assume that $\gamma\neq 0$ 
and  consider $x, y\in (0, \infty)$ with $|y-x|<(x\wedge y)/4$. Then
\begin{equation}\label{e:B3-6.2}
|\sB(x,y)-\sB(x,x)|=\left| \left(\frac{y}{x}\right)^{\gamma}-1\right|=
\frac{|\gamma| u^{\gamma-1}|y-x|}{x^{\gamma}},
\end{equation}
where $u$ is  between $x$ and $y$.
If $x<y$, then $x=x\wedge y$ and $y\le (5/4)x$. If $\gamma\ge 1$, then the right-hand side above is less than $(5/4)^{\gamma-1}\gamma (y-x)/x$. If $\gamma <1$, then we estimate the right-hand side with 
$|\gamma| (y-x)/x$. 
Thus in both cases \textbf{(B3)} holds with $\theta=1$. If $y<x$, then we replace $x^{\gamma}$ with $y^{\gamma}$ in \eqref{e:B3-6.2} and argue as before.

See Proposition \ref{p:expectation-non-symmetry} for more on this example.

\medskip
Given an arbitrary $\sB(x,y)$ satisfying \textbf{(B1)}-\textbf{(B3)},  and the jump kernel $j(x,y)$ of the censored process $X$,  
we now construct a a pssMp $\overline{X}=(\overline{X}_t, \P_x)$ of index $\alpha$ corresponding to the jump kernel $J(x,y)=\sB(x,y)j(x,y)$ via the Lamperti transform of a certain  L\'evy process. 
Define
$$
\mu^{\sB}(y):=\sB(1,e^y)\mu(y)=\sB(1,e^y)\left(c_+\frac{e^y}{(e^y-1)^{1+\alpha}}\1_{(y>0)}+c_-\frac{e^y}{(1-e^y)^{1+\alpha}}\1_{(y<0)}\right). 
$$
By the assumptions
\textbf{(B2)}(a) and (c) we have that 
$\int_{|y|>1}\mu^{\sB}(y)dy <\infty$, while by 
\textbf{(B2)}(b) and (c) we get that 
$\int_{|y|\le 1}
y^2\mu^{\sB}(y)dy <\infty$. Thus $\mu^{\sB}$ is a L\'evy measure. 
Further, let  $\overline{\xi}$ denote the L\'evy process with infinitesimal generator 
\begin{equation}\label{e:inf-gen-overline-xi}
\overline{\AA} f(x)= -\overline{b}f'(x)+\int_{\R}\left(f(x+y)-f(x)-f'(x)y \1_{[-1,1]}(y)\right) \mu^{\sB}(y)dy,
\end{equation}
where  $\overline{b}\in \R$. 
Let $\overline{X}=(\overline{X}_t, \P_x)$ be the pssMp of index $\alpha$ obtained from $\overline{\xi}$ through the Lamperti transform. 
By using a calculation similar to the one we used to obtain $\LL^*$ in Section \ref{s:prelim}, together with the homogeneity of 
$\sB$ and \textbf{(B1)}, 
we see that the infinitesimal generator of $\overline{X}$ is 
\begin{align}\label{e:inf-gen-overline-X}
\overline{\LL} f(x)&= -\overline{b} x^{1-\alpha}f'(x)+x^{-\alpha}\int_{\R}\left(f(xe^y)-f(x)-xf'(x)y\1_{[-1,1]}(y)\right)
\mu^{\sB}(y)dy\\
&=-\overline{b} x^{1-\alpha}f'(x)+
\int_0^\infty
\left(f(z)-f(x)-xf'(x)(\log z/x)\1_{[-1,1]}(\log(z/x))\right)\sB(x,z)j(x,z)dz.\nonumber
\end{align}
This shows that the jump kernel of $\overline{X}$ is precisely $J(x,y)=\sB(x,y)j(x,y)$,
see the last sentence of Subsection \ref{ss:sspcp}.

\subsection{Symmetric case}\label{ss:2-symmetric}
In this subsection we assume that $\eta$ is a \emph{symmetric}  $\alpha$-stable process. Then $\rho=1/2$,  $c_+=c_-=:c$, and $a=0$. For simplicity, we will assume that $c=1$. 
We first consider the case that 
$\sB(x,y)$ is identically 1.
Recall that $X^*$ is the process $\eta$ killed upon exiting $(0,{\infty})$, and 
the constant $b$ in the linear term of its infinitesimal generator in  \eqref{e:linear-term}  equal to
$$
b= -
\int_0^\infty
\left((\log u)\1_{[-1,1]}(\log u)-(u-1)\1_{[-1,1]}(u-1)\right)|u-1|^{-1-\alpha}\, du.
$$

\begin{lemma}\label{l:lin-term-alt}
It holds that
$$
-b=\lim_{\epsilon\to 0} \int_{\R, |e^y-1|>\epsilon} y \1_{[-1,1]}(y)\mu(y) dy =\mathrm{p.v.} \int_{-1}^1 y \mu(y) dy.
$$
\end{lemma}
\pf We first note that by using symmetry, for $\epsilon\in (0,1)$ we have
$$
\int_{(0,{\infty}), |u-1|>\epsilon}  (u-1) \1_{[-1,1]}(u-1)|u-1|^{-1-\alpha} du=\int_{\R}v \1_{(\epsilon<|v|\le 1)}|v|^{-1-\alpha}dv=0.
$$
Therefore
\begin{align*}
I(\epsilon):= & \int_{\R, |e^y-1|>\epsilon} y \1_{[-1,1]}\mu(y) dy= 
\int_{\R, |e^y-1|>\epsilon} 
y \1_{[-1,1]}(y)\frac{e^y}{|e^y-1|^{1+\alpha}} dy\\
 =& 
 \int_{(0, \infty), |u-1|>\epsilon} 
 (\log u) \1_{[-1,1]}(\log u)|u-1|^{-1-\alpha}du\\
= & 
\int_{(0, \infty), |u-1|>\epsilon} 
\left(  (\log u) \1_{[-1,1]}(\log u)-(u-1)\1_{[-1,1]}(u-1)\right)|u-1|^{-1-\alpha}du.
\end{align*}
By letting $\epsilon \to 0$ we obtain that $\lim_{\epsilon\to 0}I(\epsilon)=-b$ which is the first equality in the statement. For the second, 
\begin{eqnarray*}
I(\epsilon)&=&
\int_{-1}^{\log(1-\epsilon)} 
y\mu(y)\, dy+\int_{\log(1+\epsilon)}^1 y\mu(y)\, dy\\
&=&\left(\int_{-1}^{\log(1-\epsilon)} y\mu(y)\, dy+\int_{-\log(1-\epsilon)}^1 y\mu(y)\, dy\right)+ \int_{\log(1+\epsilon)}^{-\log(1-\epsilon)}y\mu(y)\, dy\\
&=:&I_1(\epsilon)+I_2(\epsilon).
\end{eqnarray*}
Suppose 
$\alpha\in [1,2)$ (for $\alpha\in (0,1)$ the integral $I_2(\epsilon)$ is convergent). For $y\in (0,1/2)$ it holds that $y\mu(y)\le c_1 y^{-\alpha}$ for some $c_1>0$, hence 
\begin{align*}
I_2(\epsilon) &\le c_1\int_{\log(1+\epsilon)}^{-\log(1-\epsilon)} y^{-\alpha}dy\\
&\le c_2 
\begin{cases}
\left(\log(1+\epsilon)^{1-\alpha}-(-\log(1-\epsilon))^{1-\alpha}\right)\le c_3 \epsilon^{2-\alpha}& \text{for } \alpha\in (1,2)\\
\log\left(\frac{-\log(1-\epsilon)} {\log(1+\epsilon)} \right)
&\text{for } \alpha=1
\end{cases}
\to 0
\end{align*}
as $\epsilon \to 0$. 
Since we have already proved that
$\lim_{\epsilon\to 0}I(\epsilon)$ exists
we can conclude that
$$
\lim_{\epsilon\to 0}I_1(\epsilon)=\mathrm{p.v.} \int_{-1}^1 y \mu(y) dy.
$$
\qed

\begin{remark}\label{r:lin-term-alt}{\rm
The existence of the principal value integral $\mathrm{p.v.} \int_{-1}^1 y \mu(y) dy$ can be alternatively proved in the following way. First note that
$$
\int_{\R,\epsilon < |y|\le 1}y \mu(y) dy=\int_{\epsilon}^1 y(\mu(y)-\mu(-y))dy.
$$
Secondly, $\mu(y)-\mu(-y)=y^{-\alpha}((1-\alpha)+O(y^2))$ as $y\downarrow 0$, showing that the right-hand side above is convergent.
}
\end{remark}

Let $\overline{\xi}$ be a L\'evy process with L\'evy density $\mu^{\sB}(y)=\mu(y)\sB(1,e^y)$ and linear term 
\begin{equation}\label{e:linear-term-sym}
\overline{b}=b-\int_{-1}^1y(\sB(1,e^y)-1)\mu(y)dy,
\end{equation}
cf.~\eqref{e:inf-gen-overline-xi}. 
Note that the integral is convergent because of 
\textbf{(B3)}. 
Let $\overline{X}$ be the  corresponding pssMp  of index $\alpha$.
The jump kernel of $\overline{X}$ is $J(x,y)=\sB(x,y)|x-y|^{-1-\alpha}$ 
and the infinitesimal generator of $\overline{X}$ is given in \eqref{e:inf-gen-overline-X}. 
The following is an analog of Lemma \ref{l:lin-term-alt}.

\begin{lemma}\label{l:lin-term-alt-2}
It holds that
$$
-\overline{b}=\lim_{\epsilon\to 0} \int_{\R, |e^y-1|>\epsilon} y \1_{[-1,1]}(y)\mu^{\sB}(y) dy =\mathrm{p.v.} \int_{-1}^1 y \mu^{\sB}(y) dy.
$$
\end{lemma}
\pf 
We have
\begin{align*}
& \int_{\R, |e^y-1|>\epsilon} y \1_{[-1,1]}\mu^{\sB}(y) dy= 
\int_{\R, |e^y-1|>\epsilon} 
y \1_{[-1,1]}(y)\sB(1,e^y)\frac{e^y}{|e^y-1|^{1+\alpha}} dy\\
& = 
\int_{\R, |e^y-1|>\epsilon} y \1_{[-1,1]}(y)\sB(1,1)\frac{e^y}{|e^y-1|^{1+\alpha}} dy\\
& \quad +\int_{\R, |e^y-1|>\epsilon} y \1_{[-1,1]}(y)(\sB(1,e^y)-\sB(1,1))\frac{e^y}{|e^y-1|^{1+\alpha}} dy\\
& = :J_1(\epsilon)+J_2(\epsilon).
\end{align*}
By Lemma \ref{l:lin-term-alt}, and since $\sB(1,1)=1$, $\lim_{\epsilon\to 0}  J_1 (\epsilon)= -b$. On the other hand, by using \textbf{(B3)} if $\alpha \ge 1$, we conclude that 
$$
\lim_{\epsilon\to 0}J_2(\epsilon)
=\int_{\R} y\1_{[-1,1]}(y)(\sB(1,e^y)-1)\frac{e^y}{|e^y-1|^{1+\alpha}}dy = \int_{-1}^1 y(\sB(1,e^y)-1)\mu(y)dy.
$$
This proves the first equality in the statement.  For the second statement, note that for 
 $u\mapsto \sB(1,u)$ is by \textbf{(B3)} bounded in a neighborhood of 1. Hence, $y\mu^{\sB}(y)=y\mu(y)\sB(1,e^y)\le c_1 y^{-\alpha}$, and we obtain the conclusion in the same way as in Lemma \ref{l:lin-term-alt}. 
In the case $\alpha\in (0,1)$, since the integral is absolutely convergent, we use the dominated convergence theorem.
\qed

In the context of pssMps it is natural to write the generator in the form \eqref{e:inf-gen-overline-X} which involves a cutoff function. 
On the other hand, in the multidimensional setting of regional non-local operators,
such as the infinitesimal  
generator of a censored $\alpha$-stable process, generators are usually written as principal value integrals. 
In the context of jump kernels decaying at the boundary, such operators were studied in \cite[Section 3.2]{KSV21} when $\sB$ is symmetric (see \textbf{(B4)} below). 
In the next result
we reconcile these two approaches in the current setting.
Let 
$$
\tilde{\LL}f(x):=\textrm{ p.v. }
\int^\infty_0
(f(z)-f(x))J(x,z)dz=\lim_{\epsilon\to 0} 
\int_{(0, \infty), |z-x|>\epsilon}
(f(z)-f(x))\sB(x,z)|x-z|^{-1-\alpha} dz.
$$

\begin{lemma}\label{l:LLL}
If $f \in C_c^2((0,{\infty}))$, then $\tilde{\LL}f(x)$ is well defined and $\tilde{\LL}f=\overline{\LL}f$.
\end{lemma}
\pf
By  \textbf{(B2)}(c),
for any compact set 
$K\subset (0, \infty)$ and $\epsilon>0$,
$$\int_{z \in K, |z-x|\ge \epsilon}\sB(x,z)dz \le c(x, K, \epsilon)<\infty.$$
Using this and \textbf{(B3)},
one can follow the proofs of \cite[Lemma 3.3 and Proposition 3.4]{KSV21} and show that
$\tilde{\LL}f$ is well defined for $f \in C_c^2((0,{\infty}))$.

By the change of variables $z=xe^y$ we have:
\begin{align*}
&
\int_{(0, \infty), |z-x|>\epsilon}
(f(z)-f(x))\sB(x,z)|x-z|^{-1-\alpha} dz\\
&=\int_{\R, |xe^y-x|>\epsilon}\left(f(xe^y)-f(x)\right)\sB(x, xe^y)|x-xe^y|^{-1-\alpha}xe^y\, dy\\
&=x^{-\alpha} \int_{\R, |e^y-1|>\epsilon/x}\left(f(xe^y)-f(x)\right)\sB(1, e^y)|1-e^y|^{-1-\alpha}e^y\, dy\\
&=x^{-\alpha} \int_{\R, |e^y-1|>\epsilon/x}\left(f(xe^y)-f(x)\right)\mu^{\sB}(y)\, dy\\
&=x^{-\alpha}\left(\int_{\R, |e^y-1|>\epsilon/x}\Big(f(xe^y)-f(x)-xf'(x)y\1_{[-1,1]}(y)\Big)\mu^{\sB}(y)\, dy \right.\\
& \qquad \qquad \qquad \left. +xf'(x)\int_{\R, |e^y-1|>\epsilon/x}y\1_{[-1,1]}(y)\mu^{\sB}(y)\, dy\right)\\
&=: x^{-\alpha}(J_1(\epsilon)+xf'(x) J_2(\epsilon)).
\end{align*}
By the dominated convergence  theorem,
$$
\lim_{\epsilon \to 0}J_1(\epsilon) =\int_{\R}\Big(f(xe^y)-f(x)-xf'(x)y\1_{[-1,1]}(y)\Big)\mu^{\sB}(y)\, dy.
$$
Since $\tilde{\LL} f(x)$ is well defined, we see that there also exists
$$
\lim_{\epsilon\to 0}J_2(\epsilon)=\lim_{\epsilon\to 0}\int_{\R, |e^y-1|>\epsilon/x}y\1_{[-1,1]}(y)\mu^{\sB}(y)\, dy=: -
\tilde{b}.
$$
Thus
$$
\tilde{\LL} f(x)=-
\tilde{b}
 x^{1-\alpha} f'(x)+ x^{-\alpha}\int_{\R}\Big(f(xe^y)-f(x)-xf'(x)y\1_{[-1,1]}(y)\Big)\mu^{\sB}(y)\, dy.
$$
By Lemma \ref{l:lin-term-alt-2} we see that $
\tilde{b}=\overline{b}$ and thus $\tilde{\LL}=\overline{\LL}$. \qed

\medskip
Now we turn to the question of the behavior of the pssMp $\overline{X}$ at its 
absorption time. We assume that \begin{equation}\label{e:B-cond-exp}
\int_{-\infty}^{-1}|y|e^y \sB(1,e^y)dy + \int_1^{\infty}ye^{-\alpha y}\sB(1,e^y)dy <\infty.
\end{equation}
Then $\int_{\R, |y|\ge 1}|y|\mu^{\sB}(y)dy<\infty$, hence $\overline{\xi}_1$ has finite expectation given by
\begin{equation}\label{e:B-exp}
\E \overline{\xi}_1=-\overline{b}+\int_{\R, |y|\ge 1}y\mu^{\sB}(y)dy,
\end{equation}
cf.~\cite[Theorem 25.3, Example 25.12]{Sat14}.

For $\gamma\in \R$ let
$$
\sigma_{\gamma}(x):=\frac{e^{(1+\gamma)x}}{(e^x-1)^{1+\alpha}}-\frac{e^{-(1+\gamma)x}}{(1-e^{-x})^{1+\alpha}}
=\frac{e^{-x} (e^{(\gamma-\alpha+1)x}-1)}{(1-e^{-x})^{1+\alpha}}
, \quad x>0.
$$
The next lemma follows immediately  from the second expression of $\sigma_{\gamma}$ above.
\begin{lemma}\label{l:sign-of-sigma}
For every $x>0$ it holds that $\sigma_{\gamma}(x)>0$ for $\alpha < 1+\gamma$, $\sigma_{\gamma}(x)=0$ for $\alpha=1+\gamma$, and $\sigma_{\gamma}(x)>0$ for $\alpha >1+\gamma$.
\end{lemma}

In the next result,  we will also assume that, in addition to \textbf{(B1)}-\textbf{(B3)}, $\sB$ satisfies 

\medskip
\noindent
\textbf{(B4)} Symmetry: $\sB(x,y)=\sB(y,x)$ for all $x,y>0$.

\begin{prop}\label{p:expectation-symmetry}
Let $\overline{X}$ be a pssMp with the infinitesimal generator $\overline{\LL}$ given in \eqref{e:inf-gen-overline-X} where the jump kernel is  $\sB(x,y)|x-y|^{-1-\alpha}$ and the linear term given in \eqref{e:linear-term-sym}.  Assume that $\sB$ satisfies \textbf{(B1)}-\textbf{(B4)} and \eqref{e:B-cond-exp}. 
Let $\overline\xi$ be the corresponding L\'evy process through the Lamperti transform. Then ${\bf E} \overline\xi_1>0$ if $\alpha\in (0,1)$, ${\bf E} \overline\xi_1=0$ if $\alpha=1$, and ${\bf E} \overline\xi_1<0$ if $\alpha\in (1,2)$.
\end{prop}
\pf 
Note that by Lemma \ref{l:lin-term-alt-2} and the fact that $\sB(1,e^y)=\sB(e^y,1)=\sB(1, e^{-y})$, it holds that
$$
-\overline{b}=\mathrm{p.v.}\int_{-1}^1 y\mu^{\sB}(y)dy=\lim_{\epsilon\to 0}\int_{\epsilon}^1y(\mu^{\sB}(y)-\mu^{\sB}(-y))dy= \lim_{\epsilon\to 0}\int_{\epsilon}^1 y\sigma_0(y)\sB(1,e^y)dy.
$$
Similarly,
$$
\int_{\R, |y|\ge 1}y
\mu^{\sB}(y)dy
=\int_1^{\infty}y (\mu^{\sB}(y)-\mu^{\sB}(-y))dy =\int_1^{\infty}y \sigma_0(y)\sB(1,e^y)dy.
$$
The claim now follows from \eqref{e:B-exp} and  Lemma \ref{l:sign-of-sigma}. 
\qed

 We can also cover some cases with non-symmetric $\sB(x,y)$. 

\begin{prop}\label{p:expectation-non-symmetry}
Let $\overline{X}$ be a pssMp with the infinitesimal generator $\overline{\LL}$ given in \eqref{e:inf-gen-overline-X} where the jump kernel   $J(x,y)   =(y/x)^{\gamma}|x-y|^{-1-\alpha}$  with $\gamma\in (-1,\alpha)$ and the linear term  is  given in \eqref{e:linear-term-sym}.   
Let $\overline\xi$ be the corresponding L\'evy process through the Lamperti transform. Then $ {\bf E}  \overline\xi_1>0$ if $\alpha\in (0,1+\gamma)$, ${\bf E}  \overline\xi_1=0$ if $\alpha=1+\gamma$, and ${\bf E} \overline\xi_1<0$ if $\alpha\in (1+\gamma,2)$.

In particular, if $\gamma\in [\alpha/2,\alpha)$ then $ {\bf E} \overline\xi_1>0$.
\end{prop}
\pf 
Since $\sB(x,y)=(y/x)^{\gamma}$ with $\gamma\in (-1, \alpha)$, we have that \eqref{e:B-cond-exp} holds and that
$\mu^{\sB}(y)-\mu^{\sB}(-y)=e^{\gamma y}\mu(y)-e^{-\gamma y}\mu(-y)=\sigma_\gamma(y)$ for $y >0$. By Lemma \ref{l:lin-term-alt-2}, 
it holds that
$$
-\overline{b}=\mathrm{p.v.}\int_{-1}^1 y\mu^{\sB}(y)dy=\lim_{\epsilon\to 0}\int_{\epsilon}^1y(\mu^{\sB}(y)-\mu^{\sB}(-y))dy= \lim_{\epsilon\to 0}\int_{\epsilon}^1 y\sigma_\gamma (y)dy.
$$
Similarly,
$$
\int_{
\R, |y|\ge 1}y
\mu^{\sB}(y)dy
=\int_1^{\infty}y (\mu^{\sB}(y)-\mu^{\sB}(-y))dy =\int_1^{\infty}y \sigma_\gamma(y)dy.
$$
The claim now follows from Lemma \ref{l:sign-of-sigma}. \qed 
 
We end this subsection with a class of examples of modifying functions, 
 satisfying \textbf{(B1)}-\textbf{(B4)}, which appeared in our papers  \cite{KSV21, KSV22}
on the potential theory of Dirichlet forms with jump kernels decaying at the boundary.
For $\beta\ge 0$ and $\gamma\ge 0$ with $\gamma=0$ if $\beta=0$, we 
define
\begin{align}\label{wt_B}
\wt{B} (x,y)=\left(\frac{x\wedge y}{x\vee y}\right)^{\beta}\left(\log\left(1+\frac{x\vee y}{x\wedge y}\right)\right)^{\gamma}.
\end{align}
It is easy to check that  $\wt{B} (x,y)$ satisfies \textbf{(B1)}-\textbf{(B4)}.
Since 
$ \frac{x\wedge y}{|x-y|}\wedge 1 \asymp \frac{x\wedge y}{x\vee y}$ and  $\frac{x\vee y}{|x-y|}\wedge 1\asymp 1$,
$\wt{B}(x,y)$ is comparable to the $\sB(x,y)$
 in  \cite[(1.8)]{KSV21} with 
$\beta=\beta_1$, $\gamma=\beta_3$ and $\beta_2=\beta_4=0$

When $\sB(x,y)$ is equal to $c\wt{B}(x,y)$,  we have that the L\'evy measure $\mu^{\sB}(y)$ of $\xi$ is equal to $c$ times 
\begin{eqnarray*}
&&\frac{e^y}{|e^y-1|^{1+\alpha}}\left(\frac{1\wedge e^y}{1\vee e^y}\right)^{\beta}\left(\log\left(1+\frac{1\vee e^y}{1\wedge e^y}\right)\right)^{\gamma}\\
&=&\frac{e^y}{|e^y-1|^{1+\alpha}}\left(\1_{(y<0)} e^{y\beta}\left(\log(1+e^{-y})\right)^{\gamma}+\1_{(y>0)}e^{-y\beta}\left(\log(1+e^{y})\right)^{\gamma}\right)\\
&=& \frac{e^y}{|e^y-1|^{1+\alpha}}e^{-|y|\beta}\left(\log(1+e^{|y|})\right)^{\gamma}. 
\end{eqnarray*}
If $\gamma=0$ (so there is no logarithmic term), we see that $\mu^{\sB}$ is 
the L\'evy measure of a Lamperti stable process in the sense of \cite{CPP10}.

Clearly, if 
$\sB(x,y)$ is comparable to $\wt{B}(x,y)$, 
then the corresponding L\'evy measure $\mu_{\ast}^{\sB}$ satisfies 
$$
\mu_{\ast}^{\sB}(y)\asymp \mu^{\sB}(y), \quad y\in \R\setminus \{0\}.
$$


\section{Appendix}\label{s:appendix}

\noindent
\textbf{Proof of Lemma \ref{l:q-symmetric}:} We claim that for any $0\le j\le k$,
\begin{equation}\label{e:unique2}
\int_{(0, \infty)}
u^j(1+xu)^{-1-\alpha-k}
m(du)=0, \quad x>0.
\end{equation}
\eqref{e:unique2} is valid for $k=0$ by assumption.  
Note that, for $k=1, 2, \dots$,
$$
|\frac{\partial^k}{\partial x^k}(1+xu)^{-1-\alpha}|\le (1+xu)^{-1-\alpha}, \quad x>0, u>0.
$$
Combining this with the integrability assumption of the lemma, we can 
exchange the order of the differentiation 
and integration when we take the derivative of the left hand side of \eqref{e:unique1}.
Taking derivative with respect to $x$ in \eqref{e:unique1} we get
\begin{equation}\label{e:unique-add}
\int_{(0, \infty)}
u(1+xu)^{-1-\alpha-1}
m(du)=0, \quad x>0,
\end{equation}
and so 
\eqref{e:unique2}
 is valid for $k=j=1$. Since
\begin{align*}
&
\int_{(0, \infty)}
(1+xu)^{-1-\alpha-1}
m(du)\\
&=
\int_{(0, \infty)}
(1+xu)^{-1-\alpha}m(du)-
\int_{(0, \infty)}
xu(1+xu)^{-1-\alpha-1}m(du)=0, \quad x>0,
\end{align*}
(where the last equality follows from the assumption and \eqref{e:unique-add}), 
we get that 
\eqref{e:unique2} is valid for $k=1$ and $j=0$. Now suppose that
\eqref{e:unique2} is valid for $0\le j\le k$. Taking derivative with respect to $x$ in 
\eqref{e:unique2}, we get
$$
\int_{(0, \infty)}
u^{j+1}(1+xu)^{-1-\alpha-k-1}
m(du)=0, \quad x>0.
$$
Thus \eqref{e:unique2} is valid for $1\le j\le k+1$. Noting that
\begin{align*}
&
\int_{(0, \infty)}
(1+xu)^{-1-\alpha-k-1}m(du)\\
&=
\int_{(0, \infty)}
(1+xu)^{-1-\alpha-k}
m(du)-
\int_{(0, \infty)}
xu(1+xu)^{-1-\alpha-k-1}m(du)=0, \quad x>0, 
\end{align*}
we get that \eqref{e:unique2} is valid for $0\le j\le k+1$.

Taking $x=1$, we get that for any $0\le j\le k$,
$$
\int_{(0, \infty)}
\frac{u^j}{(1+u)^k}(1+u)^{-1-\alpha}
m(du)=0.
$$
Since the linear span of the set $\{\frac{u^j}{(1+u)^k}: 0< j< k\}$
is an algebra of real-valued continuous functions on $(0, \infty)$ which separates
points of $(0, \infty)$ and vanishes at infinity, by the Stone-Weierstrass Theorem,
the linear span of the set 
$\{\frac{u^j}{(1+u)^k}: 0< j< k\}$ 
 is dense in $C_\infty(0, \infty)$ with respect to the uniform topology. Thus for all $g\in C_\infty(0, \infty)$,
$$
\int_{(0, \infty)}
g(u)(1+u)^{-1-\alpha}m(du)=0,
$$
which implies 
$(1+u)^{-1-\alpha}m(du)$ is the zero measure on $(0, \infty)$. Therefore
$m$ is the zero measure on $(0, \infty)$.
\qed

\vspace{.1in}
\textbf{Acknowledgment}:
We thank the referees for insightful comments and suggestions that led to improvements of the paper. 
We also thank Pierre Patie for helpful comments on a preliminary version of this paper.

\vskip 0.1truein

\parindent=0em

{\bf Panki Kim}

Department of Mathematical Sciences and Research Institute of Mathematics,

Seoul National University, Seoul 08826, Republic of Korea

E-mail: \texttt{pkim@snu.ac.kr}

\bigskip

{\bf Renming Song}

Department of Mathematics, University of Illinois, Urbana, IL 61801,
USA

E-mail: \texttt{rsong@math.uiuc.edu}

\bigskip

{\bf Zoran Vondra\v{c}ek}

Department of Mathematics, Faculty of Science, University of Zagreb, Zagreb, Croatia,

Email: \texttt{vondra@math.hr}

\end{document}